\documentclass[12pt,leqno]{amsart}
\usepackage{amssymb} 
\begin{document}
\theoremstyle{plain}
\newtheorem{MainThm}{Theorem}
\newtheorem{thm}{Theorem}[section]
\newtheorem{clry}[thm]{Corollary}
\newtheorem{prop}[thm]{Proposition}
\newtheorem{lem}[thm]{Lemma}
\newtheorem{deft}[thm]{Definition}
\newtheorem{hyp}{Assumption}
\newtheorem*{conjecture}{Conjecture}

\theoremstyle{definition}
\newtheorem{rem}[thm]{Remark}
\newtheorem*{acknow}{Acknowledgements}
\numberwithin{equation}{section}
\newcommand{\eps}{\varepsilon}
\renewcommand{\phi}{\varphi}
\renewcommand{\d}{\partial}
\newcommand{\re}{\mathop{\rm Re} }
\newcommand{\im}{\mathop{\rm Im}}
\newcommand{\R}{\mathbf{R}}
\newcommand{\C}{\mathbf{C}}
\newcommand{\N}{\mathbf{N}} 
\newcommand{\Z}{\mathbf{Z}} 
\newcommand{\D}{C^{\infty}_0} 
\renewcommand{\O}{\mathcal{O}}
\newcommand{\dbar}{\overline{\d}}
\newcommand{\supp}{\mathop{\rm supp}}
\newcommand{\abs}[1]{\lvert #1 \rvert}
\newcommand{\csubset}{\Subset}
\newcommand{\detg}{\lvert g \rvert}
\newcommand{\msetminus}{\setminus}

\title[Limiting Carleman weights]{Limiting Carleman weights and anisotropic inverse problems}
\author[]{David Dos Santos Ferreira, Carlos E. Kenig, Mikko Salo, Gunther Uhlmann}
\begin{abstract}
In this article we consider the anisotropic Calder\'on problem and related inverse problems. The approach is based on limiting Carleman weights, introduced in \cite{KSU} in the Euclidean case. We characterize those Riemannian manifolds which admit limiting Carleman weights, and give a complex geometrical optics construction for a class of such manifolds. This is used to prove uniqueness results for anisotropic inverse problems, via the attenuated geodesic X-ray transform. Earlier results in dimension $n \geq 3$ were restricted to real-analytic metrics.
\end{abstract}
\maketitle
\setcounter{tocdepth}{1} 
\setcounter{secnumdepth}{2}
\tableofcontents
%
%
\begin{section}{Introduction and main results}

\begin{subsection}{Introduction}

In this paper we consider the Calder\'on problem in the anisotropic case.
This inverse method, also called Electrical Impedance Tomography (EIT),
consists in determining the conductivity of a medium by making
voltage and current measurements at the boundary. Applications range from
geophysical prospection to medical imaging.

Anisotropic conductivities depend on direction. Muscle tissue in the human
body is an important example of an anisotropic conductor. For instance
cardiac muscle has a conductivity of 2.3 mho in the transverse direction
and 6.3 mho in the longitudinal direction \cite{BB}. The conductivity in this case is
represented by a positive definite, smooth, symmetric matrix $\gamma=(\gamma^{jk}(x))$ 
in a domain $\Omega$ in Euclidean space.

If there are no sources or sinks of current in $\Omega$, the potential $u$ in $\Omega$, given a voltage potential $f$ on $\partial \Omega$, solves the Dirichlet problem
\begin{align*} 
   \left\{ 
   \begin{aligned}
      \frac{\partial}{\partial x_j} \left( \gamma^{jk} \frac{\partial u}{\partial x_k} \right) &= 0 & & \textrm{ in } \Omega, \\
      u&=f & &\textrm{ on } \d \Omega.
   \end{aligned}
   \right.
\end{align*}
Here and throughout this article we are using Einstein's summation convention: repeated indices in lower and upper position are summed. The boundary measurements are given by the Dirichlet-to-Neumann map (DN map), defined by
\begin{equation*} 
\Lambda_{\gamma} f = \gamma^{jk} \frac{\partial u}{\partial x_j} \nu_k \Big|_{\partial \Omega}
\end{equation*}
where $\nu=(\nu_1, \ldots, \nu_n)$ denotes the unit outer normal to $\partial \Omega$ and $u$ is the solution of the Dirichlet problem.  The inverse problem is whether one can determine $\gamma$ by knowing $\Lambda_{\gamma}$.

Unfortunately, $\Lambda_{\gamma}$ doesn't determine $\gamma$ uniquely. This observation is due to L.~Tartar (see \cite{KV} for an account). Let $\psi:\overline{\Omega} \to \overline{\Omega}$ be a $C^\infty$ diffeomorphism with $\psi|_{\partial \Omega}=\text{Id}$ where $\text{Id}$ is the identity map. Then 
\begin{equation*}
\Lambda_{\widetilde\gamma}=\Lambda_{\gamma}
\end{equation*}
where
\begin{equation*} 
\widetilde\gamma=\left (\frac{{}^t\psi' \cdot \gamma \cdot \psi'}{|\det \psi'|}\right )\circ\psi^{-1}.
\end{equation*}
Here $\psi'$ denotes the (matrix) differential of $\psi$, ${}^t\psi'$ its transpose, and the dot $\cdot$ represents multiplication of matrices.

We have then a large number of conductivities with the same DN map: any change of variables of $\Omega$ that leaves the boundary fixed gives rise to a new conductivity with the same boundary measurements. The question is whether this is the only obstruction to unique identifiability of the conductivity. It is known that this is the case in two dimensions. The anisotropic problem can be reduced to the isotropic one by using isothermal coordinates (Sylvester \cite{Sy}), and 
combining this with the result of Nachman \cite{N} for isotropic conductivities gives the 
result for anisotropic conductivities with two derivatives. The regularity was improved by Sun and Uhlmann \cite{SuU2} to Lipschitz conductivities using the techniques of Brown and Uhlmann \cite{BrU}, and to $L^\infty$ conductivities by Astala-Lassas-P\"aiv\"arinta \cite{ALP} 
using the work of Astala-P\"aiv\"arinta \cite{AP}.

In the case of dimension $n\ge 3$, as was pointed out in \cite{LeU}, this is a problem of geometrical nature. 
In this article we will focus on the geometric problem.

Let $(M,g)$ be a compact Riemannian manifold with boundary. All manifolds will be assumed smooth (which means $C^{\infty}$) and oriented. The Laplace-Beltrami operator associated to the metric $g$ is given in local coordinates by 
   $$ \Delta_g u = \detg^{-1/2}  \frac{\partial}{\partial x_{j}} \left(\detg^{1/2} \, g^{jk} \frac{\partial u}{\partial x_{k}} \right) $$
where as usual $(g^{jk})$ is the matrix inverse of $(g_{jk})$, and $\detg = \det(g_{jk})$. Let us consider the Dirichlet problem 
\begin{equation*} 
\Delta_ g u  = 0\hbox{ in }M,\quad
u|_{\partial M} =  f.
\end{equation*}
The DN map in this case is defined as the normal derivative 
\begin{equation*} 
\Lambda_g f = \partial_{\nu} u|_{\partial M} = g^{jk} \frac{\partial u}{\partial x_j} \nu_k  \Big|_{\partial M}
\end{equation*}
where $\nu=\nu^l \d_{x_l}$ denotes the unit outer normal to $\partial M$, and $\nu_k = g_{kl} \nu^l$ is the conormal. 
The inverse problem is to recover $g$ from $\Lambda_g$.

There is a similar obstruction to uniqueness as for the conductivity. We have 
\begin{equation} \label{isometry_obstruction}
\Lambda_{\psi^\ast g} = \Lambda_g
\end{equation}
where $\psi$ is a $C^\infty$ diffeomorphism of $M$ which is the
identity on the boundary. As usual $\psi^\ast g$ denotes the pull back of
the metric $g$ by the diffeomorphism $\psi$. 

In the two dimensional case there is an additional
obstruction since the Laplace-Beltrami operator is conformally invariant.
More precisely 
\[
\Delta_{c g}=\frac{1}{c}\Delta_g
\]
for any function $c$, $c \ne 0$. Therefore we have that for $n=2$
\begin{equation} \label{twodim_obstruction}
\Lambda_{c(\psi^\ast g)}=\Lambda_{g}
\end{equation}
for any smooth function $c \ne 0$ so that $c |_{\partial M}=1$.

Lassas and Uhlmann \cite{LaU} proved that \eqref{isometry_obstruction} is the only obstruction to
unique identifiability of the metric for real-analytic manifolds in
dimension $n\ge 3$. In the two dimensional case they showed that \eqref{twodim_obstruction} is
the only obstruction to unique identifiability for smooth Riemannian surfaces.  Moreover
these results assume that $\Lambda_g$ is measured only on an open subset of the
boundary.

Notice that these two results don't assume any condition on the topology of
the manifold except for connectedness. An earlier result of Lee-Uhlmann \cite{LeU} assumed
that $(M,g)$ was strongly convex and simply connected. 
The result of \cite{LaU} in dimension $n \ge 3$ was extended by Lassas-Taylor-Uhlmann \cite{LTU} to non-compact, connected
real-analytic manifolds with boundary.
Einstein manifolds are real-analytic in the interior and it was conjectured in
\cite{LaU} that Einstein manifolds are determined, up to isometry, from the DN map.
This was recently proven by Guillarmou-Sa Barreto \cite{GS}. 

These results on the anisotropic Calder\'on problem for $n \geq 3$ are based on the analyticity of the metric. The recovery of the metric in the interior of $M$ proceeds by analytic continuation, using the knowledge of Taylor series of $g$ at the boundary. Thus, these results do not give information from the interior of the manifold.

On the other hand, in the isotropic case where $g$ is a conformal multiple of the Euclidean metric, many results are available even for nonsmooth coefficients. These results are based on special complex geometrical optics solutions to elliptic equations, introduced in \cite{SU}. These have the form 
\begin{equation*}
u = e^{-\frac{1}{h}\langle \zeta, x \rangle}(1+r_0),
\end{equation*}
where $\zeta \in \C^n$ is a complex vector satisfying $\zeta \cdot \zeta = 0$, and $r_0$ is small as $h$ tends to $0$. However, complex geometrical optics solutions have not been available in the anisotropic case, which has been a major difficulty in the study of that problem.

One of the main contributions of this paper is a complex geometrical optics construction for a class of Riemannian manifolds. 
This is based on the work of Kenig-Sj\"ostrand-Uhlmann \cite{KSU}, where more general complex geometrical optics solutions, 
of the form 
\begin{equation*}
u = e^{-\frac{1}{h}(\varphi + i\psi)}(a+r_0),
\end{equation*}
were constructed in Euclidean space. Here $\varphi$ is a limiting Carleman weight.

In this paper we characterize those Riemannian manifolds which admit limiting Carleman weights, and also characterize all such weights in Euclidean space. We give a construction of complex geometrical optics solutions on a class of Riemannian manifolds, and we use these solutions to prove uniqueness results in inverse problems. The inverse problems considered are the recovery of an electric potential and a magnetic field from boundary measurements on an admissible Riemannian manifold, and the determination of an admissible metric within a conformal class from the DN map. Let us now state the precise results of this article.

\end{subsection}

\begin{subsection}{Statement of results}

We first recall the definition of limiting Carleman weights. Let $h > 0$ be a small parameter, and consider the semiclassical Laplace-Beltrami operator $P_0=-h^2\Delta_g$. If $\varphi$ is a smooth real-valued function on $M$, consider the conjugated operator 
\begin{align}
\label{Intro:Conjugated}
   P_{0,\phi} = e^{\phi/h} P_0 e^{-\phi/h}. 
\end{align}
Here it is natural to work with open manifolds (i.e.~manifolds without boundary such that no component is compact).

\begin{deft}
   A real-valued smooth function $\phi$ in an open manifold $(M,g)$ is said to be a \emph{limiting Carleman weight} if it has non-vanishing 
   differential, and if it satisfies on $T^* M$ the Poisson bracket condition 
   \begin{align}
   \label{Intro:LCW}
      \{\overline{p_{\phi}},p_{\phi}\}=0 \ \textrm{ when } \ p_{\phi}=0,
   \end{align} 
   where $p_{\phi}$ is the principal symbol, in semiclassical Weyl quantization, of the conjugated Laplace-Beltrami
   operator~$(\ref{Intro:Conjugated})$.
\end{deft}   

Our first result is a characterization of those Riemannian manifolds which admit limiting Carleman weights.

\begin{MainThm}
\label{Intro:CharLCW}
   If $(M,g)$ is an open manifold having a limiting Carleman weight, 
   then some conformal multiple of the metric $g$ admits a parallel unit vector field.
   For simply connected manifolds, the converse is also true.
\end{MainThm}

Locally, a manifold admits a parallel unit vector field if and only if it is isometric to the product of an Euclidean interval and another Riemannian manifold. This is an instance of the de Rham decomposition \cite{Pe}, or is easy to prove directly (see Lemma \ref{appendix:localLem}). Thus, if $(M,g)$ has a limiting weight $\varphi$, one can choose local coordinates in such a way that $\phi(x)=x_1$ and
      $$ g(x_1,x') = c(x)\left( \begin{array}{cc} 1 & 0 \\ 0 & g_0(x') \end{array} \right), $$
where $c$ is a positive conformal factor. Conversely, any metric of this form admits $\varphi(x) = x_1$ as a limiting weight.

In the case $n=2$, limiting Carleman weights in $(M,g)$ are exactly the harmonic functions with non-vanishing differential (see Section 2). The case $n \geq 3$ is more complicated.
However, for the Euclidean metric it is possible to determine all the limiting Carleman weights.

\begin{MainThm}
\label{Intro:EuclLCW}
   Let $\Omega$ be an open subset of $\R^n$, $n\geq 3$, and let $e$ be the Euclidean metric.
   The limiting Carleman weights in $(\Omega,e)$ are locally of the form
      $$ \phi(x)=a\phi_0(x-x_0)+b $$
   where $a \in \R \msetminus \{0\}$ and $\phi_0$ is one of the following functions:
   \begin{align*}
      \langle x,\xi \rangle, & \quad  \arg \langle x,\omega_1 + i \omega_2 \rangle, \\  
      \log|x|, \quad \frac{\langle x,\xi \rangle}{|x|^2}, & \quad \arg \big( e^{i\theta} (x+i\xi)^2 \big), \quad
      \log \frac{|x+\xi|^2}{|x-\xi|^2}
   \end{align*}
   with $\omega_1, \omega_2$ orthogonal unit vectors, $\theta \in [0,2\pi)$ and $\xi \in \R^n \msetminus \{0\}$. 
\end{MainThm}

We use the following definition for the argument function 
    $$ \arg z = 2 \arctan \frac{\im z}{|z|+\re z}, \quad z \in \C \setminus \R_-. $$

\begin{rem}
   The possible weights are all real analytic functions. Some comments are made at the end of section \ref{Euclsection} about the global aspect
   of Theorem \ref{Intro:EuclLCW}.
\end{rem}

Let us now introduce the class of manifolds which admit limiting Carleman weights and for which we can prove uniqueness results in inverse problems. For this we need the notion of simple manifolds \cite{Sh}.
\begin{deft}
   A manifold $(M,g)$ with boundary is \emph{simple} if $\partial M$ is strictly convex%
   \footnote{cf. Definition \ref{appendix:strictlyconvex}.}
   , and for any point $x \in M$ the exponential map
   $\exp_x$ is a diffeomorphism from some closed neighborhood of $0$ in $T_x M$ onto $M$.
\end{deft}
\begin{deft}
   A compact manifold with boundary $(M,g)$, of dimension $n \geq 3$, is \emph{admissible}
   if it is conformal to a submanifold with boundary of $\R \times (M_0,g_0)$ where $(M_0,g_0)$ is a compact
   simple $(n-1)$-dimensional manifold.
\end{deft}
Examples of admissible manifolds include the following:
\begin{enumerate}

\item[1.] 
Bounded domains in Euclidean space, in the sphere minus a point, or in hyperbolic space. In the last two cases, the manifold is conformal to a domain in Euclidean space via stereographic projection. 

\item[2.] 
More generally, any domain in a locally conformally flat manifold is admissible, provided that the domain is appropriately small. Such manifolds include locally symmetric 3-dimensional spaces, which have parallel curvature tensor so their Cotton tensor vanishes (see the Appendix).

\item[3.] 
Any bounded domain $M$ in $\R^n$, endowed with a metric which in some coordinates has the form 
\begin{equation*}
   g(x_1,x') = c(x) \left( \begin{array}{cc} 1 & 0 \\ 0 & g_0(x') \end{array} \right),
\end{equation*}
with $c > 0$ and $g_0$ simple, is admissible.

\item[4.] 
The class of admissible metrics is stable under $C^2$-small perturbations of $g_0$.

\end{enumerate}

The first inverse problem involves the Schr\"odinger operator 
\begin{equation*}
\mathcal{L}_{g,q} = -\Delta_g + q,
\end{equation*}
where $q$ is a smooth complex valued function on $(M,g)$. We make the standing assumption that $0$ is not a Dirichlet eigenvalue of $\mathcal{L}_{g,q}$ in $M$. Then the Dirichlet problem 
\begin{align*} 
   \left\{ 
   \begin{aligned}
      \mathcal{L}_{g,q} u &= 0 & & \textrm{ in } M, \\
      u&=f & &\textrm{ on } \d M
   \end{aligned}
   \right.
\end{align*}
has a unique solution for any $f \in H^{1/2}(\partial M)$, and we may define the DN map 
\begin{equation*}
\Lambda_{g,q}: f \mapsto \partial_{\nu} u|_{\partial M}.
\end{equation*}
Given a fixed admissible metric, one can determine the potential $q$ from boundary measurements.

\begin{MainThm}
\label{Intro:MainThm}
    Let $(M,g)$ be admissible, and let $q_1$ and $q_2$ be two smooth functions on $M$. 
    If $\Lambda_{g,q_1}=\Lambda_{g,q_2}$, then $q_1=q_2$.
\end{MainThm}

This result was known previously in dimensions $n\ge 3$ for the Euclidean metric \cite{SU}
and for the hyperbolic metric \cite{I}. We remark that in the two dimensional case global
uniqueness is not known even for the Euclidean metric. It is known for potentials coming from conductivities \cite{N} or for a generic class of potentials \cite{SuU}. 

We obtain similar uniqueness results for the Schr\"odinger operator in the presence of a magnetic field.
Let $A$ be a smooth complex valued 1-form on $M$ (the magnetic potential), and denote
\begin{align*}
   \mathcal{L}_{g,A,q} &= {d_{\bar{A}}}^*d_A+q,
\end{align*}
where $d_A=d+iA \wedge \; : C^{\infty}(M) \to \Omega^1(M)$ and ${d_A}^*$ is the formal adjoint of 
$d_A$ (for the sesquilinear inner product induced by the Hodge dual on the exterior form algebra). This reads in local coordinates
\begin{align*}
   \mathcal{L}_{g,A,q}u = -\detg^{-1/2} \big(\d_{x_{j}}+iA_{j}\big) \big(\detg^{1/2} g^{jk} \big(\d_{x_{k}}+iA_k\big)u\big) 
\end{align*}
if $A=A_j \,dx^j$.

As before, we assume throughout that $0$ is not a Dirichlet eigenvalue of $\mathcal{L}_{g,A,q}$ in $M$, and consider the Dirichlet problem 
\begin{align*}
   \left\{ 
   \begin{aligned}
      \mathcal{L}_{g,A,q}u &= 0 & & \textrm{ in } M, \\
      u&=f & &\textrm{ on } \d M.
   \end{aligned}
   \right.
\end{align*} 
We can define the DN map as the magnetic normal derivative 
\begin{equation*}
\Lambda_{g,A,q}: f \mapsto d_A u(\nu)|_{\partial M}.
\end{equation*}

This map is invariant under gauge transformations of the magnetic potential: we have 
      $$ \Lambda_{g,A+d\psi,q} = \Lambda_{g,A,q} $$
for any smooth function $\psi$ which vanishes on the boundary. Thus, it is natural that one recovers the magnetic field $dA$ and electric potential $q$ from the map $\Lambda_{g,A,q}$.

\begin{MainThm}
\label{Intro:ThmMagn}
    Let $(M,g)$ be admissible, let $A_1,A_2$ be two smooth 1-forms on $M$ 
    and let $q_1,q_2$ be two smooth functions on $M$. If $\Lambda_{g,A_1,q_1}=\Lambda_{g,A_2,q_2}$, then
    $dA_1=dA_2$ and $q_1=q_2$.
\end{MainThm}

This result was proved in \cite{NSU} for the Euclidean metric. Our proof is closer to \cite{DSFKSU} which considers partial boundary measurements. See \cite{Sa} for further references on the inverse problem for the magnetic Schr\"odinger operator in the Euclidean case.

The next result considers the anisotropic Calder\'on problem. Under the additional condition that the metrics are in the same conformal class, one expects uniqueness since the only diffeomorphism that leaves a conformal class invariant is the identity. In dimensions $n \geq 3$ this was known earlier for metrics conformal to the Euclidean metric \cite{SU}, conformal to the hyperbolic metric \cite{I}, and analytic metrics in the same conformal class \cite{Li} (based on \cite{LeU}).

\begin{MainThm}
\label{Intro:ConfThm}
   Let $(M,g_1)$ and $(M,g_2)$ be two admissible Riemannian manifolds in the same conformal class.
   If $\Lambda_{g_1}=\Lambda_{g_2}$, then $g_1=g_2$.
\end{MainThm}

This article is organized as follows. In the next two sections, we study limiting Carleman weights and prove in particular
Theorems \ref{Intro:CharLCW} and \ref{Intro:EuclLCW}. In Sections 4 and 5, we prove Carleman estimates and construct complex geometrical optics solutions to the Schr\"odinger equation.
Section 6 deals with the proofs of Theorems \ref{Intro:MainThm}, \ref{Intro:ThmMagn}, and \ref{Intro:ConfThm}. 
The last two sections are devoted to two results needed in the resolution of the anisotropic inverse problems. The first is the injectivity of an attenuated geodesic X-ray transform on simple manifolds, and the second states that the DN map determines the Taylor expansion of the different quantities involved at the boundary. Finally, there is an appendix containing basic definitions and facts in Riemannian geometry which are used in this article.

\end{subsection}

\subsection*{Acknowledgements}

C.E.K.~is partly supported by NSF grant DMS-0456583. M.S.~is supported in part by the Academy of Finland. G.U. would like to acknowledge partial support of NSF and a Walker Family Endowed Professorship.  We would like to 
express our deepest thanks to Johannes Sj\"ostrand who made substantial 
contributions to this paper. His unpublished notes on characterizing 
limiting Carleman weights in the Euclidean case are the basis for sections 
2 and 3. In particular he proved that the level sets of limiting Carleman 
weights in the Euclidean case are either hyperspheres or hyperplanes (see 
section 3). We would also like to thank David Jerison for helpful discussions on limiting Carleman weights, and Robin Graham for 
useful suggestions on conformal geometry.
\end{section}
%
%
\begin{section}{Limiting Carleman weights}

We refer the reader to the appendix for a short overview on Riemannian geometry.  
We use $\langle \cdot,\cdot \rangle$ and $|\cdot|$ to denote the Riemannian inner product and norm both on the tangent 
and the cotangent space and $D$ to denote the Levi-Civita connection. Throughout this paper semiclassical conventions are 
used; we refer to \cite{DS} for an exposition of this theory.
The principal symbol of the conjugated semiclassical Laplace-Beltrami operator (\ref{Intro:Conjugated}) is given by
\begin{align}
   p_{\phi}=|\xi|^{2}-|d\phi|^{2}+2i\langle \xi,d\phi \rangle.
\end{align}

There are two reasons to use limiting Carleman weights in the construction of complex geometrical optics solutions 
   $$ u = e^{-\frac{1}{h}(\phi+i\psi)} (a+r_0) $$
of the Schr\"odinger equation $(-\Delta_g + q) u=0$. Given a function $\phi$, the construction amounts to looking for solutions of the conjugated equation
   $$ P_{0,\phi}v+h^2qv=0 $$
of the form $v=e^{-\frac{i}{h}\psi} (a+r_0)$ and then applying the usual WKB method. This includes solving the eikonal equation
   $$ p_\phi(x,d\psi)=0 $$  
and a transport equation on $a$. Note that $P_{0,\phi}$ is not a self-adjoint operator and that the symbol
$p_{\phi}$ is complex valued. The existence of a solution $\psi$ to the eikonal equation implies
   $$ \{\overline{p_{\phi}},p_{\phi}\}(x,d\psi)=0. $$
Hence using limiting Carleman weights is a way to ensure that the former (necessary) equality
is fulfilled. The other reason lies in the fact that one wants the conjugated operator $P_{0,\phi}$ to be 
locally solvable in the semiclassical sense, in order to find the remainder term $r_0$ and go from an approximate solution to an exact solution.
This means the principal symbol $p_{\phi}$ of the conjugated operator needs to satisfy H\"ormander's local solvability condition
   $$ \{\overline{p_{\phi}},p_{\phi}\} \leq 0 \textrm{ when } \ p_{\phi}=0. $$
Since applications of the complex geometrical optics construction to inverse problems require to construct solutions with both exponential 
weights $e^{\phi/h}$ and $e^{-\phi/h}$ (in order to cancel possible exponential behaviour in the product of two solutions
--- see section \ref{InvPb}) and since $p_{-\phi}=\overline{p_{\phi}}$, it seems natural to impose the bracket condition (\ref{Intro:LCW}).

We now proceed to the analysis of limiting Carleman weights on open manifolds. 
We first observe that the notion of limiting Carleman weight relates to a conformal class of Riemannian manifolds.  
\begin{lem}
\label{LCW:ConfInvLem}
   Let $(M,g)$ be an open Riemannian manifold, and $\phi$ a limiting Carleman weight.
   If $c$ is a smooth positive function, then $\phi$ is a limiting Carleman weight in $(M,c g)$.
   In particular, if $f:(\tilde{M},\tilde{g}) \to (M,g)$ is a conformal transformation, then $f^*\phi$ is a limiting Carleman 
   weight in $(\tilde{M},\tilde{g})$.
\end{lem}
\begin{proof}
   The claim follows from the fact that the principal symbol of the conjugated Laplace-Beltrami operator $e^{\phi/h}\Delta_{cg}e^{-\phi/h}$
   is $\tilde{p}_{\phi}=c^{-1} p_{\phi}$ and that
   \begin{align*}
      \frac{1}{2i} \{\overline{\tilde{p}_{\phi}},\tilde{p}_{\phi}\} =  \frac{c^{-2}}{2i}\{\overline{p_{\phi}},p_{\phi}\}
      +  c^{-1} \im \big(p_{\varphi} \{\overline{p_{\phi}},c^{-1}\}\big).
   \end{align*}
   Both term in the right-hand side vanish when $p_{\phi}=0$.
\end{proof}
\begin{rem}
   This lemma gives a way to construct limiting Carleman weights.
   If we already know a limiting Carleman weight $\phi$, then any function of the form $\phi \circ f$, where $f$ is a conformal 
   transformation on $(M,g)$, is a limiting Carleman weight. 
   
   In particular, from the linear Carleman weight $\langle x, \xi\rangle$,
   using the inversion $x \to x/|x|^2$ which is a conformal transformation on $\R^n \msetminus \{0\}$ 
   (endowed with the Euclidean metric $e$), 
   we obtain another limiting Carleman weight $\langle x,\xi \rangle/|x|^2$ on $M=\R^n \msetminus \{0\}$. 
\end{rem}

Now we know that the existence of limiting Carleman weights only depends on a conformal class of geometries.
Let $a$ and $b$ denote respectively the real and imaginary parts of the 
principal symbol of the conjugated operator (\ref{Intro:Conjugated}), so
\begin{align*}
   p_{\phi}=a+ib.
\end{align*}
They are given by the following expressions%
\footnote{The musical notation is recalled in the appendix.}
\begin{align}
\label{LCW:ab}
   a&=|\xi|^{2}-|d\phi|^{2}=|\xi^{\sharp}|^2-|\nabla \phi|^2, \\ \nonumber
   b&=2\langle d\phi,\xi \rangle =2 \langle \nabla \phi,\xi^{\sharp} \rangle.
\end{align}
Note that $\{\overline{p_{\phi}},p_{\phi}\}/2i=\{a,b\}$ so that the condition (\ref{Intro:LCW}) also reads 
\begin{align}
\label{LCW:abCond}
   \{a,b\}=0 \quad \textrm{when } a=b=0.
\end{align}

We start with the computation of this Poisson bracket in terms of the Hessian $D^2\phi$.
\begin{lem}
\label{LCW:BracketLem}
The Poisson bracket may be expressed as
\begin{align*}
   \{a,b\}(x,\xi)= 4D^2\phi(\xi^{\sharp},\xi^{\sharp}) + 4D^2\phi(\nabla \phi,\nabla \phi).
\end{align*}
\end{lem}
\begin{proof}
   Consider $a_1=|\xi|^{2}$ and $a_2=|d \phi|^2$, so that $a=a_1-a_2$, the Poisson bracket is given by 
   \begin{align*}
      \{a,b\}=\{a_{1},b\}+\{b,a_{2}\}=H_{a_1}b+H_ba_2.
   \end{align*}
   Since $a_2$ is a function only depending on $x$, and $b$ is linear in $\xi$, the second bracket
   is easily calculated and by \eqref{appendix:HessEq} we obtain
   \begin{align*}
      H_ba_2=2 L_{\nabla \phi} \big( |\nabla \phi|^2 \big) =4D^2\phi(\nabla \phi,\nabla \phi). 
   \end{align*}
   It remains to compute the first Poisson bracket. The Hamiltonian flow generated by $\frac{1}{2}a_1$ is the cogeodesic flow,
   \begin{align*}
      t \to (x(t),\xi(t)), \quad x(t) \textrm{ a geodesic}, \; \xi(t)=\dot{x}^{\flat}(t),
   \end{align*}
   therefore using \eqref{appendix:HessGeodesic} one has
   \begin{align*}
      H_{a_1}b &= 2 \frac{\d}{\d t} b\big(x(t),\xi(t) \big) \big|_{t=0} = 4 \frac{\d}{\d t} d\phi (x(t))\big(\dot{x}(t)\big) \big|_{t=0} \\
      &= 4 \frac{\d^2}{\d t^2} \phi (x(t)) \big|_{t=0} = 4D^2 \phi(\dot{x}(0),\dot{x}(0)).
   \end{align*}
   This finishes the proof of the equality since $\dot{x}(0)=\xi^{\sharp}$.
\end{proof}

The bracket condition (\ref{LCW:abCond}) now reads
\begin{multline}
\label{LCW:HessCond}
   D^2\phi(X,X)+D^2\phi(\nabla \phi,\nabla \phi)=0  \\ \textrm{ for all } X \in T(M) \textrm{ such that } |X|^2=|\nabla \phi|^2, 
   \, \langle X,\nabla \phi \rangle=0.
\end{multline} 
The situation in dimension $2$ is particularly simple.
\begin{lem}
   In the case of a $2$-dimensional Riemannian manifold $(M,g)$, limiting Carleman weights are exactly the harmonic functions
   with non-vanishing differential.
\end{lem}
\begin{proof}
   This comes from the fact that when $M$ has dimension 2, we have 
   \begin{align*}
       D^2\phi(X,X) + D^2\phi(\nabla \phi,\nabla \phi) = |\nabla \phi|^2 \mathop{\rm Tr} \, D^2 \phi
       =|\nabla \phi|^2 \, \Delta \phi
   \end{align*}
   if $|X|^2=|\nabla \phi|^2$ and $\langle \nabla \phi,X \rangle=0$.
\end{proof}

We will therefore continue our investigation in dimension $n \geq 3$. 
The expression of the Poisson bracket $\{a,b\}$ suggests that it is convenient to work with Carleman weights which are also 
distance functions, in the sense that $|\nabla \phi|=1$, since in that case $D^2\phi(\nabla \phi,\nabla \phi)=0$. One can 
always reduce to this case by using the conformal metric
\begin{align}
\label{LCW:Metric}
   \tilde{g}=|\nabla \phi|^2 g 
\end{align}
since the notion of limiting Carleman weights only depends on a conformal class of metrics.

\begin{lem}
\label{LCW:CharLem}
   Among distance functions on an open Riemannian manifold $(M,g)$, limiting Carleman weights can be characterized by the following equivalent
   properties. 
   \begin{enumerate}
      \item The Hessian of $\phi$ vanishes identically.
      \item The gradient of $\phi$ is a Killing field.
      \item The gradient of $\phi$ is a parallel field.
      \item If $x \in M$ and $v$ is in the domain of $\exp_x$ then 
              $$\phi(\exp_x v)=\phi(x) + \langle \nabla \phi(x),v\rangle.$$
   \end{enumerate}
\end{lem}
\begin{proof}
   If we assume $|\nabla \phi|=1$, then $\langle D_Z\nabla \phi,\nabla \phi\rangle=0$ for any vector field $Z$ and the bracket condition (\ref{LCW:HessCond}) is equivalent to
       $$ D^2\phi(X,X)=\langle D_X\nabla \phi,X \rangle=0 \quad \textrm{ when } \langle X,\nabla \phi\rangle =0. $$
   By bilinearity and by the fact that  $D^{2}\phi(\nabla \phi,Z) = \langle D_Z\nabla \phi,\nabla \phi\rangle=0$ one can actually drop the orthogonality 
   condition, and after polarization \eqref{LCW:HessCond} is furthermore equivalent to $D^2 \phi=0$. 
   The equivalence of (1) and (2) comes from the first equality in (\ref{appendix:HessEq}).
   The equivalence of (1) and (3) follows from the second equality in (\ref{appendix:HessEq}).
   Finally, properties (1) and (4) are equivalent because of the identities
    \begin{align*}
       \frac{\d }{\d t} \phi\big(\exp_x (tv)\big)\big|_{t=0}&=d\phi(v)=\langle \nabla \phi,v\rangle, \\
       \frac{\d^2 }{\d t^2} \phi\big(\exp_x (tv)\big)\big|_{t=0} &= D^2\phi(v,v)
    \end{align*}
    and of Taylor's expansion.
\end{proof}
\begin{rem}
   According to (4), in Riemannian manifolds, functions with null Hessian are the analogue of linear Carleman weights
   in the Euclidean setting. 
\end{rem}
\begin{rem}
\label{LCW:abRem}
   According to Lemmas \ref{LCW:BracketLem} and \ref{LCW:CharLem}, when $\phi$ is both a limiting Carleman weight
   and a distance function, the bracket $\{a,b\}$ in \eqref{LCW:abCond} vanishes everywhere instead of just on the 
   characteristic variety $a=b=0$.
\end{rem}

\begin{proof}[Proof of Theorem \ref{Intro:CharLCW}]
   Suppose that the manifold $(M,g)$ has a limiting Carleman weight~$\phi$, then $\phi$ is both a distance function and a limiting
   Carleman weight in $(M,\tilde{g})$ where $\tilde{g}$ is the metric (\ref{LCW:Metric}) conformal to $g$. According to Lemma
   \ref{LCW:CharLem}, this means that $\nabla_{\tilde{g}} \phi$ is a unit parallel field.
   
   Conversely, if $(M,g)$ is a simply connected Riemannian manifold such that $(M,c g)$ has a unit parallel field $X$, then according to Lemma
   \ref{appendix:KillLem}, $X$ is both a gradient field $X=\nabla_{c g} \phi$ and a Killing field. Thanks to Lemmas \ref{LCW:CharLem} and 
   \ref{LCW:ConfInvLem}, this implies that $\phi$ is a limiting Carleman weight in $(M,c g)$ and in $(M,g)$.
\end{proof} 

\begin{rem} \label{LCW:localRem}
It is now easy to justify the local coordinate expression in the introduction for metrics which admit limiting Carleman weights. If $(M,g)$ admits a limiting weight $\varphi$, the proof of Theorem 1 shows that $\nabla_{\tilde{g}} \varphi$ is a unit parallel field for the conformal metric $\tilde{g} = c^{-1} g$ where $c = \abs{\nabla \varphi}^{-2}$. Lemma \ref{appendix:localLem} implies that near any point of $M$ there exist local coordinates such that 
\begin{align*}
  \tilde{g}(x_1,x') = \left( \begin{matrix} 1 & 0 \\ 0 & g_0(x') \end{matrix} \right) 
  \quad \textrm{ and } \quad \nabla_{\tilde{g}} \varphi = \partial/\partial x_1 = \nabla_{\tilde{g}}(x_1).
\end{align*}
One obtains, after a translation of coordinates if necessary, that $\varphi(x) = x_1$ and 
\begin{equation*}
g(x_1,x') = c(x) \left( \begin{array}{cc} 1 & 0 \\ 0 & g_0(x') \end{array} \right).
\end{equation*}
Conversely, if $g$ is of this form and $\tilde{g} = c^{-1} g$, then $\partial/\partial x_1 = \nabla_{\tilde{g}}(x_1)$ is a unit parallel field by Lemma \ref{appendix:localLem}. Thus $\varphi(x) = x_1$ is a limiting Carleman weight.
\end{rem}

The purpose of the next lemma is to give several properties of limiting Carleman weights, in particular the fact that their level sets are totally umbilical hypersurfaces%
\footnote{cf. Definition \ref{appendix:umbilical}.}. 

\begin{lem}
\label{LCW:LemGenLCW}
   A function $\phi$ with non-vanishing differential is a limiting Carleman weight if and only if 
   $|\nabla \phi|^{-2}\nabla \phi$ is a conformal Killing field. In particular, if $\phi$ is a limiting Carleman weight then
   \begin{enumerate}
       \item the Hessian of $\phi$ is determined by the knowledge of $\nabla \phi$ and $D_{\nabla \phi}\nabla \phi$, that is 
                $$  D^2\phi= \lambda g + |\nabla \phi|^{-2} \big(d\phi \otimes (D_{\nabla \phi}\nabla\phi)^{\flat}
		    + (D_{\nabla \phi}\nabla\phi)^{\flat} \otimes d\phi\big) $$
             where $\lambda=-|\nabla \phi|^{-2} D^2\phi(\nabla \phi,\nabla \phi)$,
       \item the trace of the Hessian is given by
                 \begin{align*}
   	             \Delta \phi = \mathop{\rm Tr} D^2 \phi=(n-2)\lambda, 
  	        \end{align*}
  	   \item the level sets of $\phi$ in $(M,g)$ are totally umbilical submanifolds with normal $|\nabla \phi|^{-1}\nabla \phi$,
  	   with principal curvatures equal to 
                     $$\mu=-|\nabla \phi|^{-3}D^2\phi(\nabla \phi,\nabla \phi), $$        
       \item the eigenvalues of the Hessian are $\lambda$, $\kappa$ and $-\kappa$ with
                    $$ \kappa=\frac{|D_{\nabla \phi}\nabla \phi|}{|\nabla \phi|}. $$
   \end{enumerate}
\end{lem}
\begin{proof}
   Let $\tilde{g}$ denote the metric \eqref{LCW:Metric} conformal to $g$.
   Suppose that $\phi$ is a limiting Carleman weight in $(M,g)$ then according to Lemma \ref{LCW:CharLem}
       $$ \nabla_{\tilde{g}} \phi = |\nabla_g \phi|^{-2}\nabla_g \phi $$
   is a Killing field in $(M,\tilde{g})$. This implies that it is a conformal Killing field in $(M,g)$. 
   Conversely, if $\nabla_{\tilde{g}} \phi$ is a conformal Killing field in $(M,g)$, it is both a conformal Killing field
   and a unit gradient in $(M,\tilde{g})$. Evaluating the equality $L_{\nabla_{\tilde{g}} \phi}\tilde{g} = \gamma \tilde{g}$  
   at $(\nabla_{\tilde{g}} \phi,\nabla_{\tilde{g}} \phi)$ gives $\gamma =0$ since by~\eqref{appendix:HessEq}
     $$ L_{\nabla_{\tilde{g}}\phi}\tilde{g}(\nabla_{\tilde{g}} \phi,\nabla_{\tilde{g}} \phi)=
        2\langle \tilde{D}_{\nabla_{\tilde{g}} \phi}\nabla_{\tilde{g}} \phi,\nabla_{\tilde{g}} \phi\rangle
        = L_{\nabla_{\tilde{g}} \phi} |\nabla_{\tilde{g}} \phi|^2 = 0. $$
   This implies that $\nabla_{\tilde{g}} \phi$ is a Killing field in $(M,\tilde{g})$
   hence that $\phi$ is a limiting Carleman weight in $(M,\tilde{g})$ thanks to Lemma \ref{LCW:CharLem}, and 
   in $(M,g)$ thanks to Lemma \ref{LCW:ConfInvLem}.

   Formulas \eqref{appendix:LieForm} and \eqref{appendix:HessEq} give 
   \begin{align*}
      L_{|\nabla \phi|^{-2} \nabla \phi}g &=2|\nabla \phi|^{-2} \Big( D^2\phi - d\phi \otimes
      \frac{(D_{\nabla \phi}\nabla \phi)^{\flat}}{|\nabla \phi|^{2}} -
      \frac{(D_{\nabla \phi}\nabla \phi)^{\flat}}{|\nabla \phi|^{2}} \otimes d\phi \Big). 
   \end{align*}
   If $\phi$ is a limiting Carleman weight, then $|\nabla \phi|^{-2}\nabla \phi$ is a conformal Killing field,
   and the former expression equals to $\gamma g$. It is easy to compute $\gamma$ by evaluating the former expression
   at $(X,X)$ such that $\langle \nabla \phi,X  \rangle=0$ and $|X|^2=1$, and by using \eqref{LCW:HessCond}
      $$ \gamma =  2 |\nabla \phi|^{-2} D^2\phi(X,X) = -2|\nabla \phi|^{-4} D^2\phi(\nabla \phi,\nabla \phi). $$
   The combination of $L_{|\nabla \phi|^{-2} \nabla \phi}g = \gamma g$ and those two relations gives (1). 
   Besides, one also has
   \begin{align*}
      \gamma &= \frac{2}{n} \mathop{\rm div}\big(|\nabla \phi|^{-2}\nabla \phi\big) \\
      &=\frac{2}{n} \big(|\nabla \phi|^{-2}\Delta \phi-2|\nabla \phi|^{-4} D^2\phi(\nabla \phi,\nabla \phi)\big)
   \end{align*}
   and the two expressions of $\gamma$ yield the trace of the Hessian as in (2).
   
   The level sets of $\phi$ are submanifolds of $M$, with unit normal given by $\nu = |\nabla \phi|^{-1}\nabla \phi$.
   The tangent vectors of $\phi^{-1}(t)$ are orthogonal to $\nabla \phi$ and the second fundamental form is given by
    \begin{align*}
      \ell(X,Y)&= \langle D_X \nu,Y \rangle =|\nabla  \phi|^{-1} D^2\phi(X,Y) = \mu \langle X,Y\rangle 
   \end{align*}
   for $X$ and $Y$ satisfying $\langle X,\nabla \phi\rangle=\langle Y,\nabla \phi\rangle=0$.
   Thus the principal curvatures of $\phi^{-1}(t)$ are all equal to $\mu$.        
   
   A consequence of the expression of the Hessian in (1) is
   \begin{align}
   \label{LCW:CovGradPhi}
      D_Z\nabla \phi = \lambda Z + |\nabla \phi|^{-2}\big( \langle D_{\nabla \phi}\nabla \phi ,Z\rangle \nabla \phi 
      + \langle \nabla \phi,Z \rangle D_{\nabla \phi}\nabla \phi\big).
   \end{align}
   Since $D^{2}\phi=\lambda g$ on $\ker d\phi \cap \ker (D_{\nabla \phi}\nabla \phi)^{\flat}$, 
   $\lambda$ is an eigenvalue of the Hessian of multiplicity at least $n-2$.
   Because of the trace, the two remaining eigenvalues are opposite. Considering the orthogonal vectors
   \begin{align*}
       u = \frac{\nabla \phi}{|\nabla \phi|}+\frac{D_{\nabla \phi}\nabla \phi}{|D_{\nabla \phi}\nabla \phi|} \quad 
       v = \frac{\nabla \phi}{|\nabla \phi|} -\frac{D_{\nabla \phi}\nabla \phi}{|D_{\nabla \phi}\nabla \phi|}
   \end{align*}
   and using (\ref{LCW:CovGradPhi}), one has :
   \begin{align*}
        D_{u}\nabla \phi = \kappa u, \quad D_{v}\nabla \phi = -\kappa v.
   \end{align*}
   One of the two vectors $u,v$ is nonzero, therefore $\kappa$ and $-\kappa$ are the two remaining eigenvalues. Note that when $\nabla \phi$
   and $D_{\nabla \phi} \nabla \phi$ are dependent then $\kappa =|\lambda|$. In this case the eigenvalues are $\lambda$ (with multiplicity $n-1$) and
   $-\lambda$.
\end{proof}

\begin{rem}
It is an accepted fact in geometry that generic manifolds in dimension $n \geq 3$ should not admit nontrivial conformal Killing fields, thus would not admit limiting Carleman weights. This would be in striking contrast with the 2D case, where every metric is locally conformally flat and admits infinitely many limiting weights (the harmonic functions).
\end{rem}

In the Euclidean case, the log weight $\phi=\log |x|$ plays an interesting role (cf. \cite{KSU}, \cite{DSFKSU} and \cite{KnSa}).
The former computation allows us to give a partial answer to the question: if $\rho>0$ is a distance function (in the sense that $|\nabla \rho|=1$), 
when is $\phi=\log \rho$ a limiting Carleman weight? Indeed, one has
\begin{align*}
   D^2 \phi = \frac{D^2\rho}{\rho} - \frac{d \rho \otimes d\rho}{\rho^2}, \quad
   D_{\nabla \phi}\nabla \phi = -\frac{\nabla \rho}{\rho^3}, \quad \textrm{and} \quad \lambda=\frac{1}{\rho^2},
\end{align*}
thus $\phi=\log \rho$ is a limiting Carleman weight if and only if the metric $g$ has the form
\begin{align*}
   g = \rho D^2 \rho + d\rho \otimes d\rho.
\end{align*}
   
If $\varphi = \log \rho$ is a limiting Carleman weight, then by Remark \ref{LCW:localRem} there are coordinates near any point of $M$ so that $\log \rho = x_1$, $\abs{\nabla \varphi}^{-2} = e^{2 x_1}$, and 
\begin{align*}
   g(x) = e^{2 x_1} \left( \begin{array}{cc} 1 & 0 \\ 0 & g_0(x') \end{array} \right) = e^{2x_1} \tilde{g}(x').
\end{align*}
Conversely, any metric of this form admits a logarithmic Carleman weight $\varphi(x) = x_1 = \log \rho$ with $\rho = e^{x_1}$.
Therefore if $\log \rho$ is a limiting Carleman weight,  one computes the Christoffel symbols in the coordinates $(x_1,x')$ to be
   $$ \Gamma_{11}^1 = 1 \quad \textrm{ and } \quad \Gamma_{ij}^1 = -\tilde{g}_{ij}(x') \textrm{ if } i \neq 1 \textrm{ or } j \neq 1. $$ 
Thus if we denote $D_{b}=D_{\d_{b}}$, the curvature tensor for $b \neq 1$ and $c \neq 1$ satisfies 
\begin{align*}
  R_{1bc1} &= \langle D_1 D_b \partial_c - D_b D_1 \partial_c, \partial_1 \rangle \\
 	&= \langle D_1 (\Gamma_{bc}^d \partial_d) - D_b (\Gamma_{1c}^d \partial_d), \partial_1 \rangle \\
 	&= (\partial_1 \Gamma_{bc}^1 - \partial_b \Gamma_{1c}^1 + \Gamma_{bc}^d \Gamma_{1d}^1 - \Gamma_{1c}^d \Gamma_{bd}^1) g_{11} \\
 	&= (\Gamma_{bc}^1 - \Gamma_{1c}^d \Gamma_{bd}^1) g_{11}.
\end{align*}
A direct computation using the special form of $g$ shows that the last expression vanishes. Then the sectional curvature 
  $$ K(X,Y) = \frac{R(X,Y,Y,X)}{|X|^2 |Y|^2 - \langle X,Y \rangle^2} $$ 
vanishes whenever $X = \partial_1$. This implies that any manifold with nonvanishing sectional curvature at a point $p$ cannot admit a limiting Carleman weight 
of the form $\log \rho$ near $p$. In particular, the sphere or hyperbolic space do not admit limiting Carleman weights of this form.
\end{section}
%
%
\begin{section}{The Euclidean case}
\label{Euclsection}

So far, we know at least three examples of limiting Carleman weights defined on open subsets of the Euclidean space: the linear weight, 
the logarithmic weight, and the hyperbolic weight
   $$ \langle x,\xi \rangle, \quad \log |x|, \quad \frac{\langle x,\xi \rangle}{|x|^2}. $$
Note that the first is defined on $\R^n$ while the two others are on $\R^n \msetminus \{0\}$. The purpose of this section is
to determine all possible limiting Carleman weights for open subsets of the Euclidean space.

The following result is well known.
\begin{lem}
   The only connected totally umbilical hypersurfaces in the Euclidean space of dimension $n \geq 3$ are parts 
   of either hyperplanes or hyperspheres.
\end{lem}
\begin{proof}
   Let $\Sigma$ be a connected totally umbilical hypersurface in an open subset $\Omega \subset \R^n$, let $\nu$ denote its
   unit exterior normal, and let $\lambda$ be the common value of the principal curvatures. First, let us prove that $\lambda$ is constant along $\Sigma$. We have
   \begin{align*}
      D_X \nu = \lambda X
   \end{align*}  
   for all vector fields $X$ tangent to $\Sigma$. Therefore we deduce
   \begin{align*}
      D_XD_Y \nu - D_YD_X \nu = (L_X \lambda)Y-(L_Y \lambda)X+ \underbrace{\lambda [X,Y]}_{=D_{[X,Y]}\nu}.  
   \end{align*}
   Since the Euclidean space is flat, we have
   \begin{align*}
      D_XD_Y \nu - D_YD_X \nu - D_{[X,Y]}\nu = R(X,Y)\nu = 0      
   \end{align*}
   therefore we deduce $(L_X \lambda) Y - (L_Y \lambda) X=0$
   for all vectors tangent to $\Sigma$, which means that $\lambda$ is constant along $\Sigma$.
   
   We now consider the vector field $V=\sum_{j=1}^n x_j \d_{x_j}$ and we have 
      $$ D_X(\nu-\lambda V)= \lambda X - \lambda D_X V = 0 $$
   for all vector fields $X$ tangent to $\Sigma$. This means that $\nu -\lambda x$ is constant
   along the hypersurface. If $\lambda=0$, the normal is constant along the hypersurface and $\Sigma$ is part of a hyperplane,
   and if $\lambda \neq 0$ then $\alpha =\lambda^{-1}(\lambda x - \nu)$ is constant along the hypersurface,
   and $\Sigma$ is a part of the hypersphere $|x-\alpha|=1/|\lambda|$.
\end{proof}

This provides the additional information that the level sets of a limiting Carleman weight are parts of either
hyperspheres or hyperplanes. Let $\Omega$ be a bounded open connected subset of $\R^n$, $n \geq 3$, 
let $\phi$ be a limiting Carleman weight on $\Omega$, and consider $J=\phi(\Omega)$ which is an open interval. 
The curvature 
   $$ \mu(t) = -|\nabla \phi|^{-3} \langle \phi'' \nabla \phi,\nabla \phi \rangle \big|_{\phi^{-1}(t)} $$
of the level sets, is a smooth function of $t \in J$. 

We will begin by proving that $\phi$ is locally of one of the forms in Theorem \ref{Intro:EuclLCW}, 
thus we will suppose that $\Omega$ is a small open neighborhood of a point $x_0$ and $J$ is a small interval. 
First suppose that $\mu$ is identically $0$, i.e.~that \textit{all} level sets of $\phi$ are
hyperplanes. Possibly after doing a rotation, we may suppose that
\begin{align}
\label{Eucl:GradPhiInit}
   \nabla \phi (x_0) =(\d_{x_1}\phi(x_0),0,\dots,0) \neq 0.
\end{align} 
If $\Omega$ is small enough, this implies that we can take $(t,x')$ as coordinates on $\Omega$, 
and that the level sets of $\phi$ are of the form
   $$ x_1 = \langle \omega'(t),x'\rangle +s(t) $$
where as usual the prime notation stands for $x'=(x_2,\dots,x_n)$. We will begin with the relation
\begin{align}
\label{Eucl:LevelPlanes}
   \phi\big(\langle \omega'(t),x'\rangle +s(t),x'\big)=t.
\end{align}
Differentiation with respect to $x'$ gives
\begin{align}
\label{Eucl:NormalPlane}
   \omega' \d_{x_1}\phi  + \d_{x'}\phi=0 \quad \textrm{ on } \phi^{-1}(t) 
\end{align}  
which expresses the fact that $\nabla \phi$ is normal to the hyperplane~$\phi^{-1}(t)$.
    
Differentiating with respect to $t$ the relation \eqref{Eucl:LevelPlanes} we get the following equations
\begin{align}
   \label{Eucl:EqOne}
   \big(\langle \dot{\omega}',x' \rangle &+ \dot{s}\big) \d_{x_1}\phi =1, \\ 
   \label{Eucl:EqTwo} \textrm{and } \quad 
   \big(\langle \dot{\omega}',x' \rangle + \dot{s}\big)^2 \d^2_{x_1}\phi&+\big(\langle \ddot{\omega}',x' \rangle + \ddot{s}\big)
   \d_{x_1}\phi =0. 
\end{align}
It remains to compute $\d^2_{x_1}\phi$ in order to obtain ordinary differential equations for $\omega'$ and $s$. 

Differentiating \eqref{Eucl:NormalPlane} with respect to $x'$ gives
\begin{align*}
   \omega_j^2 \d^2_{x_1}\phi+2 \omega_j \d^2_{x_1x_j}\phi  + \d^2_{x_j}\phi = 0, \quad 2\leq j \leq n.
\end{align*}
Remember that we have supposed $\mu=0$, this implies that $\phi$ is harmonic because of (2) in Lemma \ref{LCW:LemGenLCW}.
Summing up the relations above and using the fact that $\Delta \phi=0$, we get
\begin{align}
 \label{Eucl:FirstRel}
   ({\omega'}^2-1) \d^2_{x_1}\phi+2 \langle \omega', \d^2_{x_1x'}\phi \rangle = 0.
\end{align}
Differentiating \eqref{Eucl:NormalPlane} with respect to $t$ we also get
\begin{align*}
   \big(\omega_j \d^2_{x_1}\phi  + \d^2_{x_1x_j}\phi\big) \big(\langle \dot{\omega}',x' \rangle + \dot{s}\big) 
   +\dot{\omega}_j \d_{x_1}\phi =0, \quad 2\leq j \leq n.
\end{align*}
Multiplying the former relations by $\omega_j$ and summing up yields
\begin{align}
\label{Eucl:SndRel}
   \big({\omega'}^2 \d^2_{x_1}\phi  + \langle \omega',\d^2_{x_1x'}\phi \rangle\big) 
   \big(\langle \dot{\omega}',x' \rangle + \dot{s}\big) + \langle \dot{\omega}', \omega' \rangle \d_{x_1}\phi=0.
\end{align}
The combination of \eqref{Eucl:FirstRel} and \eqref{Eucl:SndRel} gives
\begin{align} 
\label{Eucl:EqThree}
   \big(\langle \dot{\omega}',x' \rangle + \dot{s}\big) \d^2_{x_1}\phi 
    =-\frac{2\langle \dot{\omega}',\omega' \rangle}{1+{\omega'}^2} \d_{x_1}\phi.
\end{align}
Finally using \eqref{Eucl:EqOne}, \eqref{Eucl:EqTwo} and \eqref{Eucl:EqThree}, we have
\begin{align*}
   -\frac{2\langle \dot{\omega}',\omega' \rangle}{1+{\omega'}^2}  \big(\langle \dot{\omega}',x' \rangle + \dot{s}\big)
   +\big(\langle \ddot{\omega}',x' \rangle + \ddot{s}\big)=0. 
\end{align*}
for all $t \in J$ and all $x'$ in a neighborhood of $x'_0$.

We end up with the following system of equations
\begin{align}
\label{Eucl:SystEqPlanes}
   \left\{
   \begin{aligned} 
      \ddot{\omega}' - 2 \frac{\langle \dot{\omega}',\omega'\rangle}{1+{\omega'}^2} \dot{\omega}' &=0 \\ 
      \ddot{s} - 2 \frac{\langle \dot{\omega}',\omega'\rangle}{1+{\omega'}^2} \dot{s} &=0
   \end{aligned} \right. \quad \textrm{ on } J.
\end{align}
Solving the first ordinary differential equation yields
\begin{align*}
   \dot{\omega}_j = a_j (1+{\omega'}^2)
\end{align*}
for some constant $a_j$. After a rotation in the $x'$ variable, we may suppose that $a_j=0$ for $j \geq 3$. 
Because of the form \eqref{Eucl:GradPhiInit} of the gradient of $\phi$ at $x_0$, we have
\begin{align*}
   \omega_2(t) = \tan(a t +b), \quad \omega_j=0, \quad 3 \leq j \leq n.
\end{align*}
Injecting the former solution into the equation in $s$ leads to the following differential equation
\begin{align*}
   \ddot{s} - 2a \tan(at+b) \dot{s} = 0
\end{align*}
which can be integrated: if $a=0$ then $s(t)=ct+d$ and if $a \neq 0$ then
   $$ s(t)= c \tan(at+b) +d. $$
Finally, this gives two possible types of limiting Carleman weights
\begin{align*} 
   &1. \quad a\langle x-\tilde{x}_0,\xi \rangle +b \\ 
   &2. \quad a \arg \langle x-\tilde{x}_0,\omega_1 + i \omega_2 \rangle+b
\end{align*}
with $\omega_1, \omega_2$  unit orthogonal vectors.
\begin{rem}
   The use of complex variables simplified computations in~\cite{DSFKSU}. Note that the second weight can be written 
   as 
   \begin{align*} 
      a \arg z +b, 
   \end{align*}
   with $z=\langle x-\tilde{x}_0,\omega_1\rangle+i\langle x-\tilde{x}_0,\omega_2 \rangle$ and $a,b$ real numbers.
\end{rem}

Now we assume that $\mu$ does not vanish identically on $J$, and consider a subinterval $I \subset \mu^{-1}(\R_+^*)$ 
(the study on $\mu^{-1}(\R_-^*)$ can be done by considering $-\phi$). 
When $t\in I$, the level sets $\phi^{-1}(t)$ are therefore hyperspheres of radii $r(t)=1/\mu(t)$, 
whose centers we denote by $\alpha(t)$. Both functions $r$ and $\alpha$ are smooth on $I$, and we have 
\begin{align}
\label{Eucl:level}
   \phi\big(\alpha(t)+r(t)\omega\big)=t, \quad \forall (t,\omega) \in I \times \Gamma
\end{align}
for some open subset $\Gamma$ of the unit hypersphere. The normal to the hypersphere $\phi^{-1}(t)$ is
\begin{align}
\label{Eucl:Normal}
   \omega = |\nabla \phi|^{-1} \nabla \phi\big(\alpha(t)+r(t)\omega\big).
\end{align} 
Differentiating the identity (\ref{Eucl:level}) with respect to $t$, we obtain
\begin{align*}
   \langle \nabla \phi\big(\alpha(t)+r(t)\omega\big), \dot{\alpha}+\dot{r} \omega \rangle=1
\end{align*}
which together with (\ref{Eucl:Normal}) gives
\begin{align}
\label{Eucl:Intermediate}
   \langle \nabla \phi,\omega\rangle \big(\langle \dot{\alpha},\omega \rangle +\dot{r} \big)=1.
\end{align}
Then taking $\omega=y/|y|$ in the former equality, multiplying it by $|y|^2$ and differentiating with respect to $y$ gives
\begin{align*}
   \big(\langle \dot{\alpha},\omega \rangle + \dot{r}\big) \big(\nabla \phi + r \phi'' \omega
   -\underbrace{r\langle \phi'' \omega,\omega\rangle \omega}_{=-\nabla \phi }\big) +
   \langle \nabla \phi,\omega\rangle \big(\dot{\alpha}+\dot{r}\omega \big)=2\omega
\end{align*}
and thus
\begin{align}
\label{Eucl:Hess}
   \phi'' \omega = - \frac{|\nabla \phi|^2}{r} (\dot{\alpha}+\dot{r}\omega).
\end{align}

We go back to (\ref{Eucl:Intermediate}) and differentiate it with respect to $t$
\begin{align*}
   \langle \phi''\omega,\dot{\alpha}+\dot{r}\omega \rangle \big(\langle \dot{\alpha},\omega \rangle +\dot{r} \big) + 
   \langle \nabla \phi,\ddot{\alpha}+\ddot{r}\omega\rangle = 0.
\end{align*}
Using (\ref{Eucl:Hess}) and (\ref{Eucl:Normal}), this gives
\begin{align*}
   - \frac{1}{r}\underbrace{|\nabla \phi|\big(\langle \dot{\alpha},\omega \rangle +\dot{r} \big)}_{=1} |\dot{\alpha}+\dot{r}\omega|^2
   + \langle \ddot{\alpha},\omega\rangle +\ddot{r} = 0
\end{align*}
which leads to the following system of equations
\begin{align}
   \left\{
   \begin{aligned} 
      \ddot{\alpha} - 2 \frac{\dot{r}}{r} \dot{\alpha} &=0 \\ 
      \ddot{r} - \frac{1}{r} \big(|\dot{\alpha}|^2+{\dot{r}}^2\big) &=0
   \end{aligned} \right. \quad \textrm{ on } I.
\end{align}
The first equation implies that the centers of the hyperspheres $\phi^{-1}(t)$ are moving along a line 
with fixed direction $k \in \R^n$
\begin{align}
   \alpha(t) = \alpha_0 + k \int^t r^2(s) \, ds
\end{align}
(the indefinite integral denotes a primitive on $I$). Thus we have $|\dot{\alpha}|^2=|k|^2 r^4$ and the second equation reads 
\begin{align*}
   r \ddot{r}-r^4 |k|^2-{\dot{r}}^2 =0.
\end{align*}
It is convenient to rewrite the former equation in terms of the curvature $\mu=1/r$:
\begin{align*}
   \mu \ddot{\mu} + |k|^2 - {\dot{\mu}}^2 = 0
\end{align*}
which is equivalent to 
\begin{align}
\label{Eucl:DetEq}
   \det \left( 
   \begin{array}{cc}
      \dot{\mu} & \mu \\
      \ddot{\mu} & \dot{\mu}
   \end{array} \right) = |k|^2.
\end{align}

Before proceeding to the resolution of this differential equation, let us first notice that it actually holds on $J$.
Indeed, if we differentiate it with respect to $t$, then we obtain an equation on $\mu$
\begin{align}
\label{Eucl:DDetEq}
   \det \left( 
   \begin{array}{cc}
      \ddot{\mu} & \mu \\
      \mu^{(3)} & \dot{\mu}
   \end{array} \right) = 0
\end{align} 
which holds on any subinterval $I$ of $\mu^{-1}(\R^*)$, but also evidently on the interior of $\mu^{-1}(0)$. 
This implies that the equation \eqref{Eucl:DDetEq} holds on $J$ by continuity, hence that the equation \eqref{Eucl:DetEq} 
holds on $J$ by integration.  
 
To solve the equation \eqref{Eucl:DetEq} there are two cases to consider.

\begin{subsubsection}{Case $k=0$}

This is the case where the hyperspheres $\phi^{-1}(t)$ are concentric. The vectors $(\mu,\dot{\mu})$ and its derivative are 
linearly dependent, thus (\ref{Eucl:DetEq}) is equivalent to the first order differential equation $\dot{\mu}=c \mu$, 
with $c$ constant. To be more precise, either $\mu$ vanishes identically or one can pick $t_0 \in J$ such that 
$\mu(t_0) \neq 0$, and consider the maximal interval $I$ on which $\mu$ doesn't vanish. By \eqref{Eucl:DetEq} 
the derivative of $\dot{\mu}/{\mu}$ on $I$ is zero, therefore $\mu$ solves the equation $\dot{\mu}=c \mu$ on $I$, 
with $c$ constant. Then the curvature is an exponential function of $t$ which never vanishes, thus $I=J$. 
The function $\phi$ is easily determined to be a logarithmic weight:  
\begin{align*}
   \phi(x) = a \log|x-\tilde{x}_0| + b. 
\end{align*}
\begin{rem}
\label{Eucl:Complex} 
   Note that the logarithmic weight can be written as
      $$ \phi =  a \re (\log z) +b $$ 
   with $z=y_1+i|y'|$, $y=x-x_0=(y_1,y')$.
\end{rem}
\end{subsubsection}

\begin{subsubsection}{Case $k \neq 0$}     

The equation \eqref{Eucl:DDetEq} shows that the vectors $(\mu,\dot{\mu})$ and its second derivative are 
linearly dependent, this implies that $\mu$ solves a second order differential equation of the form 
$\ddot{\mu}=c \mu$ with $c$ a constant. To see this, consider the function
   $$ \sigma = \frac{\ddot{\mu}}{\mu} $$
which is \textit{a priori} only defined outside the set of zeros of $\mu$. 
Note that if $\mu(\tau)=0$ then by \eqref{Eucl:DetEq} one has $\dot{\mu}(\tau)\neq 0$ and therefore the
zeros of $\mu$ are isolated. Moreover, by \eqref{Eucl:DDetEq} one has $\ddot{\mu}(\tau)=0$. 
We can extend $\sigma$ as a continuous function on all of $J$, because as $t$
tends to $\tau$ we have 
   $$ \sigma(t) = \frac{\ddot{\mu}(t)-\ddot{\mu}(\tau)}{t-\tau} \, \frac{t-\tau}{\mu(t)-\mu(\tau)} 
      \longrightarrow \frac{\mu^{(3)}(\tau)}{\dot{\mu}(\tau)}. $$
The Taylor expansion at higher orders actually shows that $\sigma$ is of class $C^1$.
Differentiating the function $\sigma$ on $J \setminus \mu^{-1}(0)$ gives 
   $$ \dot{\sigma} = \frac{\mu^{(3)}\mu-\ddot{\mu}\dot{\mu}}{\mu^2}=0 $$
and $\dot{\sigma}=0$ on $J$ by continuity.
   
Hence the function $\sigma$ is constant and this implies that $\mu$ satisfies the equation $\ddot{\mu}-c\mu=0$
with $c$ constant. Depending on whether $c$ is zero, negative or positive, the curvature is one of the following functions 
\begin{align}
\label{Eucl:curv}
   \pm|k|(t+b), \quad \frac{|k|}{a} \sin(at+b), \quad \pm \frac{|k|}{a}\sinh(at + b).
\end{align}
If $t_0=\phi(x_0)$ and $\mu(t_0) \neq 0$, we may choose $J$ so small that $\mu$ doesn't vanish on $J$.
If $\mu(t_0)=0$, then either $\mu$ vanishes identically near $0$, and that case was covered first,
or $t_0$ is an isolated zero of $\mu$ by \eqref{Eucl:DetEq} with $k \neq 0$, and we may choose $J$ so small that $\mu$ doesn't vanish on
$J \setminus \{t_0\}$. In conclusion, one may assume that the curvature vanishes at most once on $J$.  
The corresponding expressions for the centers of the spheres $\alpha$ are 
\begin{gather*}
   \alpha_0-\frac{k}{|k|^2} \, \frac{1}{t+b}, \quad
   \alpha_0-\frac{k}{|k|^2} \, \frac{a}{\tan(at+b)}, \\
   \alpha_0-\frac{k}{|k|^2} \, \frac{a}{\tanh(at+b)}.
\end{gather*}

If $\mu(t_0)=0$ then these expressions depend \textit{a priori} on the connected component $I_{\pm}$ 
of $J \setminus \{t_0\}$ on which $\mu$ is positive or negative. The norm of the vector $k$ however is determined by the equation 
\eqref{Eucl:DetEq} which holds on $J$ and its direction $k/|k|$ by the normal to the hyperplane $\phi^{-1}(t_0)$. 
This may be seen in the following way: consider the family of hyperspheres
    $$ \Big| x-\alpha_{\pm}+\frac{k_{\pm}}{|k|^2} \, \frac{a}{\tanh(at+b)}\Big|^2 = \frac{a^2}{|k|^2 \sinh^2(at + b)} $$
when $t \in I_{\pm}$. After expansion, we obtain
\begin{equation}
\label{Eucl:LimSpheres}
    2 \langle x-\alpha_{\pm}, k_{\pm} \rangle = \frac{|k|^2 \tanh(at+b)}{a} \Big( -\frac{a^2}{|k|^2} - |x-\alpha_{\pm}|^2 \Big)
\end{equation}
and letting $t$ tend to $t_0=-b/a$, we get the equation of the hyperplane $\phi^{-1}(t_0)$
    $$ \langle x-\alpha_{\pm}, k_{\pm} \rangle = 0. $$
This implies that either $k_+=k_-$ or $k_+=-k_-$.
Besides, by \eqref{Eucl:Normal} and \eqref{Eucl:Intermediate}, we have $\langle \dot{\alpha},\omega \rangle+\dot{r}=|\nabla \phi|^{-1}$ on $I_+$,
thus $\langle k_{+},\nabla \phi \rangle -\dot{\mu} |\nabla  \phi|= \mu^2$ and by  letting $x$ tend to $x_0$,
this gives
  $$ \langle \nabla \phi(x_0),k_{+} \rangle = \dot{\mu}(t_0)|\nabla \phi(x_0)|.$$ 
Considering $I_-$, we obtain the same identity with $k_+$ replaced by $k_-$. This additional information removes the uncertainty on the sign and we have $k_+=k_-$. 
Moreover the components of $\alpha_{+}$ and $\alpha_{-}$ along $k$ are also equal. The computations are similar 
in the two other cases.

To determine explicitly the function $\phi$, let us deal with the case where the curvature is $|k|\sin(at+b)/a$: 
to fix the ideas, suppose that $aJ+b \subset (-\theta,-\theta+2\pi)$ with $\theta \in \R$ and that 
$\mu$ vanishes at $t_0=\phi(x_0)$ (the non-vanishing case is simpler). 
Then we have $J \setminus \{t_0\} = I_- \cup I_+$ and 
   $$ x = \alpha_{\pm} - \frac{a}{|k|^2 \tan (a\phi(x)+b)} \, k + \frac{a}{|k|\sin(a \phi(x)+b)} \, \omega $$
for all $x\in \phi^{-1}(I_{\pm})$. Thus, if we take $\xi=-ak/|k|^2$, we have
\begin{align*}
   e^{i \theta} (x-\alpha_{\pm}+i\xi)^2 &= \frac{2|\xi|^2}{\sin^2 (a\phi+b)} \Big(\cos (a\phi+b)-
   \frac{\langle \omega,k \rangle}{|k|}\Big) e^{i (a\phi+b+\theta)} \\ 
   &= 2 |a|^{-1} |\xi|\big(\langle \dot{\alpha},\omega \rangle + \dot{r}\big)e^{i (a\phi+b+\theta-\pi)} 
\end{align*}
and since $a\phi+b+\theta-\pi \in (-\pi,\pi)$, and $\langle \dot{\alpha},\omega \rangle + \dot{r}>0$ by \eqref{Eucl:Normal} 
we have
   $$ a\phi(x)+b = \arg \big(e^{i \theta} (x-\alpha_{\pm}+i\xi)^2\big)+\pi-\theta $$
on $\phi^{-1}(I_{\pm})$. With the former expression of $\phi$ we have
   $$ |a\nabla \phi(x)| = \frac{2 |\xi|}{|x-\alpha_{\pm}|^2-|\xi|^2} \quad \textrm{on} 
      \quad \langle x-\alpha_{\pm},\xi \rangle =0 $$
this shows that $|x-\alpha_+|^2 = |x-\alpha_-|^2$ in a neighborhood of $x_0$ in the hyperplane $\phi^{-1}(t_0)$,
hence that $\alpha_+=\alpha_-$. The former expression of $\phi$ therefore holds on $\Omega$.
 
The two other cases are similar. This finally gives three possible types of limiting Carleman weights
\begin{align*} 
   &1. \quad a \frac{\langle x-\tilde{x}_0,\xi \rangle}{|x-\tilde{x}_0|^2}+b \\ 
   &2. \quad a \arg \big(e^{i\theta}(x-\tilde{x}_0+i\xi)^2\big)+b \\ 
   &3. \quad a \mathop{\rm arctanh} \frac{2 \langle x-\tilde{x}_0,\xi \rangle}{|x-\tilde{x}_0|^2+|\xi|^2}+b
       = a\log \frac{|x-\tilde{x}_0+\xi|^2}{|x-\tilde{x}_0-\xi|^2}+b.
\end{align*}
\begin{rem}
   As in Remark \ref{Eucl:Complex}, these functions take a simple form with respect to some complex variable.
   Take $y=x-\tilde{x}_0$, and $(y_1,y')$ such that $y=y_1 \, \xi/|\xi|+y'$ where $y'$ is orthogonal to $\xi$. 
   Denote $z=y_1+i|y'|$, the Carleman weights take the form
   \begin{align*} 
   &1. \quad a\re \frac{1}{z}+b, \\ 
   &2. \quad a \im \log \frac{e^{i\theta}(z+ic)}{z-ic}+b, \\ 
   &3. \quad a \re \log \frac{z+c}{z-c}+b, 
\end{align*}
 with $a,b,c$ real numbers.
\end{rem}
This proves Theorem \ref{Intro:EuclLCW}. Note that this result allows to determine limiting Carleman weights locally on the sphere or on 
the hyperbolic space by conformal transformation. 

\begin{rem}
\label{Intro:remLCWdomains}
   While the linear weight is smooth on $\R^n$, the weights $\log |x|$ and $\langle x,\xi\rangle/|x|^2$ are only well defined and smooth
   on $\R^n \setminus \{0\}$, and the function $\log |x+\xi|^2/|x-\xi|^2$ on $\R^n \setminus \{\xi,-\xi\}$. The weight 
   $\arg \langle x,\omega_1 + i \omega_2 \rangle$ is defined on 
   $\R^n \setminus \{\langle x,\omega_2\rangle=0, \;\langle x,\omega_1\rangle \leq 0\}$
   and the function $\arg \big(e^{i\theta}(x+i\xi)^2\big)$ on 
      $$ \R^n \setminus \Big \{ |x+\mathop{\rm cotan} \theta \, \xi|^2 = \frac{|\xi|^2}{\sin^2 \theta},
         \quad \langle x,\xi \rangle  \lesseqgtr 0 \Big\} $$
   if $\theta \in (0,\pi)$ (with the $\geq$ inequality) or $\theta \in (\pi,2\pi)$ (with the $\leq$ inequality) and on 
      $$ \R^n \setminus \big \{ \langle x,\xi \rangle = 0, \quad |x| \lesseqgtr |\xi| \big\} $$
   if $\theta = 0$ (with the $\leq$ inequality) or $\theta=\pi$ (with the $\geq$ inequality).     
\end{rem}

Let us end this section with some comments about the global aspect of Theorem \ref{Intro:EuclLCW} when $\Omega$ is an open
connected set. Let $\phi$ be a limiting Carleman weight on $(\Omega,e)$. 
The function $\phi$ is real analytic on $\Omega$: 
indeed if $x_0 \in \Omega$, the function $\phi$ is one of the limiting Carleman 
weights calculated above in the neighourhood of $x_0$, hence real analytic.
This ensures that if $\phi$ is locally equal to one of the functions
\begin{align*}
   a  \langle x-\tilde{x}_0,\xi \rangle+b, \quad &a\log|x-\tilde{x}_0|+b, \\
   a\frac{\langle x-\tilde{x}_0,\xi \rangle}{|x|^2}+b, \quad
   &a\log \frac{|x-\tilde{x}_0+\xi|^2}{|x-\tilde{x}_0-\xi|^2}+b
\end{align*}
then $\phi$ is equal to this function on the whole set $\Omega$. Indeed if $\phi$ is equal to a linear weight near 
some point $x_0 \in \Omega$ then $\phi$ is equal to this weight on the whole set $\Omega$ by analytic continuation. 
If $\phi$ is equal to the function $a \log |x-\tilde{x}_0|+b$ near 
$x_0 \in \Omega$ then $\Omega$ cannot contain the singularity $\tilde{x}_0$ of this function.
Otherwise, $\phi$ would be equal to this function on $\Omega \setminus \{\tilde{x}_0\}$ by analytic continuation 
and would not blow up at $\tilde{x}_0$.  Thus $\phi$ is equal to $a \log |x-\tilde{x}_0|+b$ on $\Omega$. 
The proof is similar if $\phi$ is one of the functions $a \langle x-\tilde{x}_0,\xi \rangle /|x-x_0|^2+b$ or
$a \log (x-\tilde{x}_0+\xi)^2/(x-\tilde{x}_0-\xi)^2+b$.
   
For the two argument forms, we need some additional assumptions on $\Omega$.
Suppose that $\phi$ is equal to $a \arg \big(e^{i\theta}(x-\tilde{x}_0+i\xi)^2\big)+b$ near $x_0$, 
and that the image of the set $\Omega$ by $x \mapsto (x-\tilde{x}_0+i\xi)^2$ is contained
in a simply connected set $U \subset \C^*$, then by analytic continuation, $\phi$ is equal to
   $$ a\arg_U (x-\tilde{x}_0+i\xi)^2+b $$ 
where $\arg_U$ is the determination of the argument on $U$ which coincides with $\arg(e^{i\theta}z)$ at 
$z_0=(x_0-\tilde{x}_0+i\xi)^2$. The point is similar for the function 
$a \arg \langle x-\tilde{x}_0,\omega_1 + i \omega_2 \rangle$.
In particular, if $\Omega$ is contained in one of the domains of existence computed in Remark \ref{Intro:remLCWdomains}
and $\phi$ is locally of the corresponding argument form, then this is still true globally on $\Omega$. 

\end{subsubsection}
\end{section}
%
%
\begin{section}{Carleman estimates}

Let $(M,g)$ be a compact Riemannian manifold with boundary. By $dV$ we denote the volume form on $(M,g)$, and by $dS=\nu \lrcorner \, dV$ the induced volume form on $\d M$.
The $L^2$ norm of a function is then given by
    $$ \|u\|_{L^2(M)}=\Big(\int_M |u|^2 \, dV\Big)^{\frac{1}{2}} $$
and the corresponding scalar product by
    $$ (u|v) = \int_M u \, \overline{v} \, dV. $$ 
Similarly on the boundary, the norm and scalar products are given by
    $$ \|f\|_{L^2(\d M)}=\Big(\int_{\d M} |f|^2 \, dS\Big)^{\frac{1}{2}} \quad (f|h)_{\d M} = \int_{\d M} f \, \overline{h} \, {dS}. $$
We write for short
    $$ \|\nabla u\|_{L^{2}(M)}= \big\| |\nabla u| \big\|_{L^2(M)} =\Big(\int_M |\nabla u|^2 \, dV\Big)^{\frac{1}{2}}$$
and we denote by $H^1_{\rm scl}(M)$ the semiclassical Sobolev space associated to the norm
    $$ \|u\|_{H^1_{\rm scl}(M)}=  \big(\|u\|_{L^2(M)}^2+ \|h \nabla u\|_{L^2(M)}^2\big)^{\frac{1}{2}}. $$
We assume that $(M,g)$ is embedded in a compact manifold $(N,g)$ without boundary, and that $\varphi$ is a limiting Carleman weight on $(U,g)$ 
where $U$ is open in $N$ and $M \Subset U$. The goal of this section is to prove the following Carleman estimate.
    
\begin{thm}
\label{Carl:CarlEstThm}
   Let $(U,g)$ be an open Riemannian manifold and $(M,g)$ a compact Riemannian submanifold with boundary such that $M \Subset U$.
   Suppose that $\phi$ is a limiting Carleman weight on $(U,g)$.
   Let $X$ be a smooth vector field on $M$ and $q$ a smooth function on $M$. There exist two constants $C>0$ and $0<h_0 \leq 1$ such that for all 
   functions $u \in \D(M^{\circ})$ and all $0<h\leq h_0$, one has the inequality 
   \begin{align}
   \label{Carl:CarlEst}
       \|e^{\frac{\phi}{h}}u\|_{H^1_{\rm scl}(M)} \leq C h \|e^{\frac{\phi}{h}}(\Delta+X+q) u\|_{L^2(M)}.
   \end{align} 
\end{thm}

To lighten the notations we will forget the subscript $L^2(M)$ whenever it is not needed.

\begin{proof}
   We first observe that the result is invariant under conformal change of metrics since we have
   \begin{align*}
      c^{\frac{n+2}{4}} \big(\Delta_{g}+X+q\big)u = (\Delta_{c^{-1}g}+c X+q_{c})\big(c^{\frac{n-2}{4}} u\big)
   \end{align*}
   with $q_{c}=c q-\frac{n-2}{4} Xc+c^{\frac{n+2}{4}}\Delta_g\big(c^{-\frac{n-2}{4}}\big)$.
   Therefore if needed, we can assume the limiting Carleman weight to be a distance function by replacing $g$ by the conformal 
   metric (\ref{LCW:Metric}).
   
   Our next observation is that the estimate may be perturbed by zero order terms since this gives rise to an error
   of the form $\O(h) \|e^{\phi/h} u\|$, which may be absorbed in the left-hand side if $h$ is assumed small enough. 
   Therefore we can neglect the potential $q$, and assume that $q=0$ from the start. Let us first assume that we also have $X=0$.
   Then the Carleman estimate (\ref{Carl:CarlEst}) is equivalent to the following \textit{a priori} estimate
   \begin{align}
   \label{Carl:Apriori}
       \|v\|_{H^1_{\rm scl}(M)} \leq C_1 h^{-1} \|P_{0,\phi}v\|.
   \end{align} 
   One goes from one inequality to another by taking $v=e^{\phi/h}u$. The conjugated operator is given by
       $$ P_{0,\phi}=-h^2\Delta-|\nabla \phi|^2+2\langle \nabla \phi,h\nabla \rangle+h\Delta \phi. $$
   Then we have in particular
   \begin{align*}
      \|h\nabla v\|^2 = ( P_{0,\phi}v|v) + \big\||\nabla \phi| \, v \big\|^2 - 2 \big(\langle \nabla \phi,h\nabla v\rangle \big| v\big) 
      - h ( \Delta \phi \, v|v)
   \end{align*}
   therefore using Cauchy-Schwarz inequality one sees that
   \begin{align}
   \label{Carl:SobNorm}
      \|h\nabla v\|^2 \leq \|P_{0,\phi}v\|^2 + C_1 \|v\|^2. 
   \end{align}
   This means that the gradient of $v$ may be  controlled and that it suffices to prove the \textit{a priori} estimate
   \begin{align*}
       \|v\| \leq C_2 h^{-1} \|P_{0,\phi}v\|
   \end{align*} 
   to obtain (\ref{Carl:Apriori}).
   
   We decompose $P_{0,\phi}$ into its self-adjoint and skew-adjoint parts
      $$ P_{0,\phi}=A+iB, \quad A=-h^2 \Delta - |\nabla \phi|^2, \quad B=2\langle\nabla \phi,h\nabla\rangle+h\Delta \phi $$
   and we have by integration by parts 
   \begin{align}
   \label{Carl:ABexp}
      \|P_{0,\phi} v\|^2 = \|Av\|^2 + \|Bv\|^2 + i ( [A,B]v|v). 
   \end{align}
   A direct application of the commutator method will not be enough to get an \textit{a priori} estimate
   assuming the bracket condition (\ref{LCW:abCond}), one needs to use convexification. This classical argument
   consists in taking a modified weight $f \circ \phi$ where $f$ is a convex function chosen
   so that the bracket in (\ref{LCW:abCond}) becomes positive. We decompose the operator
   $P_{0,f \circ \phi}=\tilde{A}+i\tilde{B}$ into its self-adjoint and skew-adjoint parts, and denote by $\tilde{a}$ and
   $\tilde{b}$ the corresponding principal symbols. We now suppose, as we may according to our first observation,
   that $\phi$ is both a limiting Carleman weight and a distance function. We have  
   \begin{align*}
      \nabla(f \circ \phi) &= (f' \circ \phi) \, \nabla \phi \\
      D^2(f \circ \phi) &= (f''\circ \phi) \, d\phi \otimes d\phi + \underbrace{(f'\circ \phi) \, D^2\phi}_{=0}
   \end{align*}
   therefore using Lemma \ref{LCW:BracketLem}
   \begin{align*}
      \big\{\tilde{a},\tilde{b}\}(x,\xi) &=  4(f'' \circ \phi) \, (f'\circ \phi)^2 |\nabla \phi|^4 
      + 4(f'' \circ \phi) \, \langle\nabla \phi,\xi^{\sharp}\rangle^2 \\ &=4(f'' \circ \phi) \, (f'\circ \phi)^2 
      + \underbrace{(f'' \circ \phi)(f'\circ \phi)^{-2}}_{=\beta} \, \tilde{b}^2.
   \end{align*}
   At the operator level, this gives
   \begin{align*}
      i[\tilde{A},\tilde{B}] = 4h (f'' \circ \phi) \, (f'\circ \phi)^2 + h \tilde{B} \beta \tilde{B} + h^2 R 
   \end{align*}
   where $R$ is a first order semiclassical differential operator. For the function $f$, we choose the following convex polynomial  
      $$ f(s) = s+\frac{h}{2\eps} s^2, \quad f'(s)=1+\frac{h}{\eps}s , \quad f''(s)=\frac{h}{\eps}. $$
   We choose $h/\eps \leq \eps_0 <1$ with $\eps_0$ small enough so that $f'> \frac{1}{2}$ on $\phi(M)$
   and denote $\phi_{\eps}=f \circ \phi$. Note that the coefficients of $R$, as well as $\beta$, are uniformly bounded with respect to $h$ 
   and $\eps$. 
   
   We finally obtain
   \begin{align*}
      i\big( [\tilde{A},\tilde{B}]v \big| v \big) &\geq \frac{h^2}{\eps} \|v\|^2 - C_3 h \|\tilde{B} v\|^2
      -C_3 h^2 \|v\|_{H^1_{\rm scl}} \|v\|
   \end{align*}
   and using (\ref{Carl:SobNorm})
   \begin{align*}
      i\big( [\tilde{A},\tilde{B}]v \big| v \big)  &\geq \frac{h^2}{\eps}(1-C_4 \eps) \|v\|^2 - C_3 h \|\tilde{B} v\|^2
      -C_3 \|P_{0,\phi_{\eps}}v\|^2.
   \end{align*} 
   Going back to (\ref{Carl:ABexp}), this gives 
   \begin{align}
   \label{Carl:X=0}
      (1+C_3) \|P_{0,\phi_{\eps}} v\|^2 &\geq \|\tilde{A}v\|^2 + (1-C_3h)\|\tilde{B}v\|^2 \\ \nonumber &\quad+ \frac{h^2}{\eps}(1-C_4 \eps) \|v\|^2. 
   \end{align}   
   If we don't assume $X=0$, the conjugated operator $P_{0,\phi_{\eps}}$ has to be perturbed by an additional term of the form
      $$ h^2 X_{\phi_{\eps}}=h^2e^{\phi_{\eps}/h}Xe^{-\phi_{\eps}/h}=h^2X-h f'\circ \phi X\phi. $$  
   By \eqref{Carl:SobNorm} and the estimate $\|h^2 X_{\phi_{\eps}} v\| \leq C_5 h\|v\|_{H^1_{\rm scl}}$, 
   the inequality (\ref{Carl:X=0}) may easily be perturbed into 
   \begin{align*}
      2(1+C_3) \|P_{0,\phi_{\eps}}v+h^2X_{\phi_{\eps}} v\|^2 \geq \frac{h^2}{\eps}(1-C_6 \eps) \|v\|^2
   \end{align*}
   if $h$ is small enough. Taking $\eps$ small enough, we obtain
   \begin{align}
   \label{Carl:Carl_tobe_shifted}
      C_{7} \|P_{0,\phi_{\eps}}v+h^2X_{\phi_{\eps}} v\|^2 \geq \frac{h^2}{\eps} \|v\|^2 
   \end{align}
   which implies with the choice of $f=s+hs^2 /2\eps$ that
      $$ C_{\eps} \|e^{\frac{\phi^2}{2 \eps}} e^{\frac{\phi}{h}} (\Delta+X) u\|^2 \geq 
          h^2 \|e^{\frac{\phi^2}{2 \eps}} e^{\frac{\phi}{h}} u\|^2 $$
   therefore we obtain the desired estimate since $1 \leq e^{\frac{\phi^2}{2 \eps}} \leq C'_{\eps}$.
\end{proof}

\begin{rem}
   The use of G{\aa}rding's inequality could give a stronger Carleman estimate, as in \cite{SaTz}.
   The present proof makes it possible to include boundary terms, which is useful in the study of
   inverse problems with partial data (see \cite{KSU}, \cite{DSFKSU} and \cite{KnSa}). 
\end{rem}

In order to prove suitable solvability results, we need to shift the indices of the Sobolev spaces in the Carleman estimate by using pseudodifferential 
calculus. Recall that we assume that $(M,g)$ is embedded in a compact manifold $(N,g)$ without boundary, and that $\varphi$ is a limiting Carleman weight 
near $(\bar{U},g)$ where $U$ is open in $N$ and $M \Subset U$.  The Laplace-Beltrami operator $-\Delta$ on $N$, with domain $C^{\infty}(N) \subset L^2(N)$, 
is essentially self-adjoint with spectrum in $[0,\infty)$. By the spectral theorem we may define for $s \in \R$ the semiclassical Bessel potentials 
\begin{equation*}
    J^s = (1-h^2 \Delta)^{s/2}.
\end{equation*}
One has $J^s J^{t} = J^{s+t}$, and $J^s$ commutes with any function of $-\Delta$. Define for $s \in \R$ the semiclassical Sobolev spaces via 
\begin{equation*}
    \|u\|_{H^s_{\text{scl}}(N)} = \|J^s u\|_{L^2(N)},
\end{equation*}
so $H^s_{\text{scl}}(N)$ is the completion of $C^{\infty}(N)$ in this norm. It is easy to see that the dual of $H^s_{\text{scl}}(N)$ may be isometrically identified with $H^{-s}_{\text{scl}}(N)$.

It is a basic fact that $J^s$ is a semiclassical pseudodifferential operator of order $s$ in $N$ (see \cite{T} and
\cite{DS}). This implies pseudolocal estimates: if $\psi, \chi \in \D(N)$ with $\chi = 1$ 
near $\supp \psi$, and if $s,\alpha,\beta \in \R$ and $K \in \N$, then 
\begin{equation}
\label{Carl:pseudolocal_est}
    \|(1-\chi)J^s \psi u\|_{H^{\alpha}_{\text{scl}}(N)} \leq C_K h^K \|u\|_{H^{\beta}_{\text{scl}}(N)}.
\end{equation}
We will also use commutator estimates in the form 
\begin{equation} 
\label{Carl:commutator_est}
    \|[A,J^s] u\|_{L^2(N)} \leq C h \|u\|_{H^{s}_{\text{scl}}(N)}
\end{equation}
whenever $A$ is a first order semiclassical differential operator in $N$.
\begin{lem}
\label{Carl:ShiftedCarl}
   Under the above assumptions on $M,N,U$ and under the assumptions of Theorem \ref{Carl:CarlEstThm}, given $s \in \R$ there are two constants $C_{s}>0$ and 
   $0<h_s \leq 1$ such that for all functions $u \in \D(M^{\circ})$ and all $0<h\leq h_s$ one has the inequality 
   \begin{align}
   \label{Carl:ShiftedCarlEst}
       \|e^{\frac{\phi}{h}}u\|_{H^{s+1}_{\rm scl}(N)} \leq C_{s} h \|e^{\frac{\phi}{h}}(\Delta_g+X+q) u\|_{H^s(N)}.
   \end{align} 
\end{lem}
\begin{proof}
    We consider the conjugated operator
         $$ P_{\phi_{\eps}}=e^{\frac{\phi_{\eps}}{h}}h^2(\Delta_g+X+q)e^{-\frac{\phi_{\eps}}{h}}$$
    where $\phi_{\eps}$ is the weight defined in the proof of Theorem \ref{Carl:CarlEstThm}. Let $\chi \in \D(U)$ with $\chi = 1$ near $M$. 
    Then the estimates \eqref{Carl:Carl_tobe_shifted}, \eqref{Carl:pseudolocal_est} imply 
    \begin{align*}
          h \|u\|_{H^{s+1}_{\text{scl}}} &\leq h \|\chi J^s u\|_{H^1_{\text{scl}}} + h \|(1-\chi) J^s u\|_{H^1_{\text{scl}}} \\
          &\leq C_1 \sqrt{\varepsilon} \|P_{\phi_{\eps}} (\chi J^s u)\|_{L^2} + C_1h^2 \|u\|_{H^{s+1}_{\text{scl}}}.
    \end{align*}
    By the estimate $\|[P_{\phi_{\eps}}, \chi] J^s u\|_{L^2} \leq C_2 h^2 \|u\|_{H^{s+1}_{\text{scl}}}$, and by
    absorbing the error terms $\|u\|_{H^{s+1}_{\text{scl}}}$ in the left hand-side if $h$ is small enough, we obtain 
    \begin{equation} 
    \label{Carl:shift_intermediate_estimate}
          h \|u\|_{H^{s+1}_{\text{scl}}} \leq C_1 \sqrt{\varepsilon} \|J^s P_{\phi_{\eps}} u\|_{L^2} + 
          C_1 \sqrt{\varepsilon} \|\chi [P_{\phi_{\eps}}, J^s] u\|_{L^2}.
    \end{equation}
    Note that in estimating the last term one may extend $\varphi$ smoothly outside $U$ if desired. Since 
    \begin{equation*}
        P_{\phi_{\eps}} = -h^2 \Delta - |\nabla \varphi_{\varepsilon}|^2 + 2 \langle \nabla \varphi_{\varepsilon}, h\nabla \rangle
        + h\Delta \varphi_{\varepsilon}+h^2X_{\phi_{\eps}} +h^2q
    \end{equation*}
    and since $[-h^2 \Delta, J^s] = 0$, the commutator estimates \eqref{Carl:commutator_est} imply that for $h \ll \varepsilon \ll 1$
    the last term in \eqref{Carl:shift_intermediate_estimate} may be absorbed in the left hand side.
    This finishes the proof since $1 \leq e^{\frac{\phi^2}{2\eps}} \leq C_{\eps}$.
\end{proof}
\begin{prop}
\label{Carl:LocalSolv}
   If $A$ is a smooth 1-form and $q$ is a smooth function on $M$, there exists a constant $0<h_0\leq 1$ such that for any function $f \in L^2(M)$
   there exists a solution $u  \in H^1(M)$ to the equation $e^{\phi/h}\mathcal{L}_{g,A,q}e^{-\phi/h}u=f$ satisfying   
   \begin{align*}
      \|u\|_{H^1_{\rm scl}(M)} \leq C h \|f\|_{L^2(M)}.
   \end{align*}    
\end{prop}
\begin{proof}
    We consider the conjugated operator
         $$ P_{\phi}^*=e^{-\frac{\phi}{h}}h^2\mathcal{L}_{g,\bar{A},\bar{q}}e^{\frac{\phi}{h}}.$$
    Let $f \in L^2(M)$, we consider the subspace $E=P_{\phi}^*(\D(N))$ of $H^{-1}_{\text{scl}}(N)$ and the linear form  defined on $E$ by
    \begin{align*}
         L(P_{\varphi}^* v) = (f|v)_{L^2(M)}, \quad v \in \D(M^{\circ}).
    \end{align*}
    By Lemma \ref{Carl:ShiftedCarlEst} applied to $P_{\varphi}^*$, this form is well defined and 
    \begin{equation*}
        |L(P_{\varphi}^* v)| \lesssim \|f\|_{L^2(M)} h^{-1} \|P_{\varphi}^* v\|_{H^{-1}_{\text{scl}}(N)}.
    \end{equation*}
    By the Hahn-Banach theorem, there is an extension $\hat{L}$ of $L$ which is a bounded functional on $H^{-1}_{\rm scl}(N)$ with norm less than $h^{-1} \|f\|_{L^2}$. 
    Since the dual of $H^{-1}_{\text{scl}}(N)$ is $H^1_{\text{scl}}(N)$, there exists a function $\tilde{u} \in H^1_{\rm scl}(N)$ such that $\hat{L}(v) = (\tilde{u}|v)$ and 
    $h \|\tilde{u}\|_{H^1_{\text{scl}}(N)} \lesssim \|f\|_{L^2(M)}$. Then $u = \tilde{u}|_M$ is the desired solution, since for all $v \in \D(M^{\circ})$ 
    \begin{equation*}
        (P_{\varphi} u|v) = (u|P_{\varphi}^* v) = \hat{L}(P_{\varphi}^* v) = L(P_{\varphi}^* v) = (f|v).
    \end{equation*}
    This completes the proof.
\end{proof}
\end{section}
%
%
%
\begin{section}{Complex geometrical optics}

Let $\varphi$ be a limiting Carleman weight in an admissible manifold $(M,g)$. We will construct solutions to $\mathcal{L}_{g,q} u = 0$ in $M$ of the form 
\begin{equation} \label{CGO:cgo_form}
u = e^{-\frac{1}{h}(\varphi+i\psi)}(a+r_0).
\end{equation}
Here the real valued phase $\psi$ and complex amplitude $a$ are obtained from a WKB construction, and the function $r_0$ will be a correction term which is small when $h$ is small.

We write $\rho = \varphi + i\psi$ for the complex phase. It will be convenient to extend the notations $\langle \,\cdot\,,\,\cdot\, \rangle$ and $\lvert \,\cdot\, \rvert^2$ to complex tangent vectors by 
\begin{gather*}
\langle \zeta, \eta \rangle = \langle \re\,\zeta, \re\,\eta \rangle - \langle \im\,\zeta, \im\,\eta \rangle + i(\langle \re\,\zeta, \im\,\eta \rangle + \langle \im\,\zeta, \re\,\eta \rangle), \\
\lvert \zeta \rvert^2 = \langle \zeta, \zeta \rangle.
\end{gather*}
We make similar extensions for the inner product of cotangent vectors. With this notation, the conjugated operator $P_{\rho} = e^{\rho/h} h^2 \mathcal{L}_{g,q} e^{-\rho/h}$ has the expression 
\begin{equation*}
P_{\rho} = -\abs{\nabla \rho}^2 + h(2(\nabla \rho) + \Delta \rho) + h^2 \mathcal{L}_{g,q}.
\end{equation*}
Then \eqref{CGO:cgo_form} will be a solution of $\mathcal{L}_{g,q} u = 0$ provided that $P_{\rho}(a+r_0) = 0$. Following the WKB method, this results in the three equations 
\begin{eqnarray}
 & \abs{\nabla \rho}^2 = 0, & \label{CGO:eikonal} \\
 & (2(\nabla \rho) + \Delta \rho) a = 0, & \label{CGO:transport} \\
 & P_{\rho} r_0 = - h^2 \mathcal{L}_{g,q} a. & \label{CGO:correction}
\end{eqnarray}
These equations will be solved in special coordinates in the admissible manifold $(M,g)$. We know that $(M,g)$ is conformally embedded in $\R \times (M_0,g_0)$ for some compact 
simple $(n-1)$-dimensional $(M_0,g_0)$. Assume, after replacing $M_0$ with a slightly larger simple manifold if necessary, that for some 
simple $(D,g_0) \csubset (\mathrm{int}\,M_0,g_0)$ one has 
\begin{equation} 
   \label{CGO:m_embedding}
    (M,g) \csubset (\R \times \mathrm{int}\,D,g) \csubset (\R \times \mathrm{int}\,M_0,g).
\end{equation}
Here $\R \times M_0$ is covered by a global coordinate chart in which $g$ has the form 
\begin{equation} 
\label{CGO:metric_global}
   g(x) = c(x) \left( \begin{array}{cc} 1 & 0 \\ 0 & g_0(x') \end{array} \right),
\end{equation}
where $c > 0$ and $g_0$ is simple. The limiting Carleman weight will be $\varphi(x) = x_1$.

\begin{subsection}{The eikonal equation}

Since $\phi$ was given, the eikonal equation \eqref{CGO:eikonal} for the complex phase becomes a pair of equations for $\psi$, 
\begin{equation*}
\abs{\nabla \psi}^2 = \abs{\nabla \phi}^2, \quad \langle \nabla \phi,\nabla \psi \rangle = 0.
\end{equation*}
One has $\varphi(x) = x_1$ and the metric is of the form \eqref{CGO:metric_global}, so $\nabla \varphi = \frac{1}{c} \frac{\partial}{\partial x_1}$ and $|\nabla \varphi| = \frac{1}{c}$. The eikonal equation now reads 
\begin{equation*}
|\nabla \psi| = \frac{1}{c}, \quad \partial_{x_1} \psi = 0.
\end{equation*}
Under the given assumptions on $(M,g)$, there is an explicit construction for $\psi$. Let $\omega \in D$ be a point such that $(x_1,\omega) \notin M$ for all $x_1$. Denote points of $M$ by $x = (x_1,r,\theta)$ where $(r,\theta)$ are polar normal coordinates in $(D,g_0)$ with center $\omega$. That is, $x' = \exp_{\omega}^D(r\theta)$ where $r > 0$ and $\theta \in S^{n-2}$. In these coordinates (which depend on the choice of $\omega$) the metric has the form 
\begin{equation*}
g(x_1,r,\theta) = c(x_1,r,\theta) \left( \begin{array}{ccc} 1 & 0 & 0 \\ 0 & 1 & 0 \\ 0 & 0 & m(r,\theta) \end{array} \right),
\end{equation*}
where $m$ is a smooth positive definite matrix.

To solve the eikonal equation it is enough to take $\psi(x) = \psi_{\omega}(x) = r$. With this choice one has $\rho = x_1 + ir$ and $\nabla \rho = \frac{2}{c} \dbar$ where 
\begin{equation*}
   \dbar = \frac{1}{2} \Big( \frac{\partial}{\partial x_1} + i \frac{\partial}{\partial r} \Big).
\end{equation*}


\end{subsection}

\begin{subsection}{The transport equation}

We now consider \eqref{CGO:transport}. In the coordinates $(x_1,r,\theta)$ one has 
\begin{equation*}
   \Delta \rho = \detg^{-1/2} \partial_{x_1} \Big( \frac{\detg^{1/2}}{c} \Big) + \detg^{-1/2} \partial_r 
   \Big(\frac{\detg^{1/2}}{c} i \Big) = \frac{1}{c} \dbar \log \frac\detg{c^2}
\end{equation*}
and the transport equation reads 
\begin{equation*}
   4 \dbar a + \Big(\dbar \log \frac\detg{c^2}\Big) a = 0.
\end{equation*}
We choose $a$ as the function 
\begin{equation*}
   a = \detg^{-1/4} c^{1/2} a_0(x_1,r) b(\theta)
\end{equation*}
where $\dbar a_0 = 0$ and $b(\theta)$ is smooth.

\end{subsection}

\begin{subsection}{Complex geometrical optics solutions}

Finally, the equation \eqref{CGO:correction} may be written as 
\begin{equation*}
e^{\varphi/h} h^2 \mathcal{L}_{g,q} e^{-\varphi/h} (e^{-i\psi/h} r_0) = - e^{-i\psi/h} h^2 \mathcal{L}_{g,q} a.
\end{equation*}
This can be solved by using Lemma \ref{Carl:LocalSolv}. We obtain $r_0$ satisfying 
   $$ \|r_0\|_{H^1_{\rm{scl}}(M)} \lesssim h. $$ 
We record the properties of the solution just obtained.

\begin{prop} 
\label{CGO:SolPotential}
   Assume that $(M,g)$ satisfies \eqref{CGO:m_embedding}, \eqref{CGO:metric_global} and let $q \in C^{\infty}(M)$. 
   Let $\omega \in D$ be such that $(x_1,\omega) \notin M$ for all $x_1$. If $(r,\theta)$ are polar normal coordinates 
   in $(D,g_0)$ with center $\omega$, then, whenever $\dbar a_0(x_1,r) = 0$ and $b(\theta)$ is smooth, the equation 
   \begin{equation*}
      \mathcal{L}_{g,q} u = 0 \ \ \text{in } M
   \end{equation*}
   has a solution of the form
   \begin{equation*}
      u = e^{-\frac{1}{h}(x_1+ir)}\big(\detg^{-1/4} c^{1/2} a_0(x_1,r) b(\theta) + r_0\big)
   \end{equation*}
   where $\|r_0\|_{H^1_{\rm{scl}}(M)} \lesssim h$.
\end{prop}

Next we consider the case where a magnetic field is present. The construction of complex geometric solutions is similar to the case 
without magnetic field, except that we have an additional factor in the amplitude.

\begin{prop} 
\label{CGO:SolMagnetic}
   Assume that $(M,g)$ satisfies \eqref{CGO:m_embedding}, \eqref{CGO:metric_global}, and let $A$ be a smooth $1$-form in $M$ and $q \in C^{\infty}(M)$. 
   Let $\omega \in D$ be such that $(x_1,\omega) \notin M$ for all $x_1$. If $(r,\theta)$ are polar normal coordinates in $(D,g_0)$ with center 
   $\omega$, then, whenever $\dbar a_0 = 0$ and $b$ is smooth, the equation 
   \begin{equation*}
      \mathcal{L}_{g,A,q} u = 0 \ \ \text{in } M
   \end{equation*}
   has a solution of the form 
   \begin{equation*}
      u = e^{-\frac{1}{h}(x_1+ir)}\big(\detg^{-1/4} c^{1/2} e^{i\Phi} a_0(x_1,r) b(\theta) + r_0\big)
   \end{equation*}
   where $\Phi$ satisfies 
   \begin{equation} 
      \label{CGO:Phi_eq}
      \dbar \Phi + \frac{1}{2}(A_1 + i A_r) = 0,
   \end{equation}
   and $\|r_0\|_{H^1_{\rm{scl}}(M)} \lesssim h$.
\end{prop}

\begin{proof}
   If $P_{\rho} = e^{\rho/h} h^2 \mathcal{L}_{g,A,q} e^{-\rho/h}$, one computes 
   \begin{equation*}
      P_{\rho} = -\lvert \nabla \rho \rvert^2 + h \big(2\nabla \rho + \Delta \rho + 2i \langle d\rho, A \rangle \big) + h^2 \mathcal{L}_{g,A,q}.
   \end{equation*}
   If $\rho = x_1+ir$ then $\lvert \nabla \rho \rvert^2 = 0$, and the transport equation for the amplitude $a$ will be 
   \begin{equation*}
      4 \dbar a + \dbar \Big(\log\,\frac{\detg}{c^2}\Big) a + 2i(A_1+iA_r) a = 0.
   \end{equation*}
   A solution is given by $a = \detg^{-1/4} c^{1/2} e^{i\Phi} a_0(x_1,r) b(\theta)$, where $\Phi$ is a solution of 
   \eqref{CGO:Phi_eq} and $\dbar a_0 = 0$. The equation for $r_0$ becomes 
   \begin{equation*}
      P_{\rho} r_0 = -h^2 \mathcal{L}_{g,A,q} a,
   \end{equation*}
   and Lemma \ref{Carl:LocalSolv} finishes the proof.
\end{proof}

\end{subsection}
\end{section}
%
%
\begin{section}{Uniqueness results}
\label{InvPb}
Let $(M,g)$ be an admissible manifold. In this section we will prove the global uniqueness results, Theorems \ref{Intro:MainThm} to \ref{Intro:ConfThm}. As often in inverse boundary problems, the starting point is an integral identity which relates the boundary measurements to solutions inside the manifold.

\begin{lem} \label{lemma:magnetic_integral_identity}
Let $A_1, A_2$ be smooth $1$-forms in $M$ and $q_1, q_2 \in C^{\infty}(M)$. If $\Lambda_{g,A_1,q_1} = \Lambda_{g,A_2,q_2}$, then 
\begin{equation} \label{magnetic_integral_identity}
   \int_{M}  \big[ i \langle A_1 - A_2, u dv - v du \rangle + \big(\abs{A_1}^2 - \abs{A_2}^2 + q_1 - q_2\big) u v \big] \,dV = 0,
\end{equation}
for any $u, v \in H^1(M)$ satisfying $\mathcal{L}_{g,A_1,q_1} u = 0$ and $\mathcal{L}_{g,-A_2,q_2} v = 0$ in $M$.
\end{lem}
\begin{proof}
One has $\Lambda_{g,A,q}^* = \Lambda_{g,\bar{A},\bar{q}}$ and 
\begin{equation*}
(\Lambda_{g,A,q} f|h)_{\partial M} = (d_A u|d_{\bar{A}} v) + (qu|v),
\end{equation*}
whenever $\mathcal{L}_{g,A,q} u = 0$ in $M$ and $u|_{\partial M} = f$, $v|_{\partial M} = h$. These facts imply the identity 
\begin{multline*}
\big((\Lambda_{g,A_1,q_1}-\Lambda_{g,A_2,q_2})f|h\big)_{\partial M} \\ = \int_{M} \big[ i \langle A_1 - A_2, u d\bar{v} - \bar{v} du \rangle
 + (\abs{A_1}^2 - \abs{A_2}^2 + q_1 - q_2) u \bar{v} \big] \,dV
\end{multline*}
where $u, v$ are $H^1(M)$ solutions of $\mathcal{L}_{g,A_1,q_1} u = 0$ and $\mathcal{L}_{g,\bar{A}_2,\bar{q}_2} v = 0$ in $M$, which satisfy $u|_{\partial M} = f$ and $v|_{\partial M} = h$. The result follows.
\end{proof}

\begin{subsection}{Recovering a potential}
%


Suppose that we have two potentials $q_1,q_2 \in C^{\infty}(M)$ for which the corresponding DN maps 
are equal. According to Lemma \ref{lemma:magnetic_integral_identity} one has
\begin{align}
\label{Aniso:ortho}
   \int_M (q_1-q_2) u_1 u_2 \, dV = 0 
\end{align}
for any $u_j \in H^1(M)$ which satisfy $\mathcal{L}_{g,q_j} u_j=0$ in $M$. We use Proposition \ref{CGO:SolPotential} and choose solutions of the form 
\begin{gather*}
   u_1 = e^{-\frac{1}{h}(x_1+ir)}\big(\detg^{-1/4} c^{1/2} e^{i\lambda(x_1+ir)} b(\theta) + r_1\big),  \\
   u_2 = e^{\frac{1}{h}(x_1+ir)}\big(\detg^{-1/4} c^{1/2} + r_2\big),  
\end{gather*}
where $\lambda$ is a real number and $\|r_j\|_{H^1_{\text{scl}}(M)} \lesssim h$. Inserting these solutions in \eqref{Aniso:ortho} 
and letting $h \to 0$ shows that 
\begin{equation*}
   \int_{\R} \iint_{M_{x_{\scriptscriptstyle 1}}} e^{i\lambda(x_1+ir)} (q_1-q_2) c(x_1,r,\theta) b(\theta) \,dr \,d\theta \,dx_1 = 0,
\end{equation*}
with $M_{x_1} = \{(r,\theta) \,;\, (x_1,r,\theta) \in M\}$. We can extend $q_1 - q_2$ smoothly by zero since $q_1 = q_2$ up to infinite order on $\partial M$ by Theorem \ref{Bdet:ThmBDet}, and we may then assume that the integral is over $\R \times D$. Taking the $x_1$-integral inside and varying $b$ gives 
\begin{equation*}
   \int e^{-\lambda r} \Big( \int_{-\infty}^{\infty} e^{i\lambda x_1} (q_1-q_2)c(x_1,r,\theta) \,dx_1 \Big) \,dr = 0, \quad \text{for all } \theta.
\end{equation*}
We denote the expression in parentheses by $f(r,\theta)$, and obtain that 
\begin{equation*}
\int_{\gamma} f(\gamma(r)) \exp\Big[-\int_0^r \lambda \,ds \Big] dr = 0
\end{equation*}
for all $D$-geodesics $\gamma$ issued from the point $\omega$. Varying $\omega$ and using the injectivity of the attenuated geodesic ray transform given in Theorem \ref{thm:attenuated_injectivity}, with constant attenuation $-\lambda$, we obtain 
\begin{equation*}
   \int_{-\infty}^{\infty} e^{i\lambda x_1} (q_1-q_2)c(x_1,r,\theta) \,dx_1 = 0, \quad \text{for all } r, \theta,
\end{equation*}
which holds for small enough $\lambda$. Since $(q_1-q_2)c$ is compactly supported in $x_1$, its Fourier transform is analytic, and we obtain that $q_1 = q_2$. This proves Theorem \ref{Intro:MainThm}.

\end{subsection}
\begin{subsection}{Recovering a magnetic field}

Next we show that also a magnetic field can be recovered from the DN map. We will need the following standard reduction of the problem to a simply connected domain.

\begin{prop} \label{prop:magnetic_extension}
   Let $M_0, M$ be compact manifolds with boundary with $M_0 \csubset M$, and let $(g,A_j,q_j)$ be smooth coefficients in $M$ such that $0$ 
   is not a Dirichlet eigenvalue of $\mathcal{L}_{g,A_j,q_j}$ in $M_0$ $(j=1,2)$. Suppose that 
   \begin{equation*}
   		A_1 = A_2, \ \ \ q_1 = q_2 \ \ \ \text{in $M \msetminus M_0$}.
   \end{equation*}
   If the DN maps $\Lambda_{g,A_j,q_j}$ in $M_0$ coincide, then the integral identity \eqref{magnetic_integral_identity} is valid 
   for any $H^1(M)$ solutions of $\mathcal{L}_{g,A_1,q_1} u = 0$ and $\mathcal{L}_{g,-A_2,q_2} v = 0$ in $M$.
\end{prop}
\begin{proof}
   Since the DN maps coincide in $M_0$, one has by Lemma \ref{lemma:magnetic_integral_identity} 
   \begin{equation*}
      \int_{M_0} \big[ i \langle A_1-A_2, u_0 dv_0 - v_0 du_0 \rangle + (|A_1|^2 - |A_2|^2 + q_1 - q_2) u_0 v_0 \,dV \big] = 0,
   \end{equation*}
   for any $H^1(M_0)$ solutions of $\mathcal{L}_{g,A_1,q_1} u_0 = 0$ and $\mathcal{L}_{g,-A_2,q_2} v_0 = 0$ in $M_0$. If $u$ and $v$ are as above, 
   then the restrictions to $M_0$ solve the corresponding equations in $M_0$, and we obtain \eqref{magnetic_integral_identity} since the coefficients 
   coincide outside $M_0$.
\end{proof}

\begin{proof}[Proof of Theorem \ref{Intro:ThmMagn}]
   By boundary determination, Theorem \ref{Bdet:ThmBDet}, and after a gauge transformation, we may extend $A_j$ and $q_j$ smoothly so that $A_1 = A_2$ and $q_1 = q_2$ 
   outside $M$. We are now in the setting of Proposition \ref{prop:magnetic_extension}. Therefore, replacing $M$ by a larger manifold inside $\R \times M_0$ if necessary, we may
   assume that $M$ is convex and $A_j$ and $q_j$ are compactly supported in $M$, and also that the identity in Lemma \ref{lemma:magnetic_integral_identity} holds whenever $u, v \in H^1(M)$ are
   solutions of $\mathcal{L}_{g,A_1,q_1} u = 0$ and $\mathcal{L}_{g,-A_2,q_2} v = 0$ in $M$.
   
   Use Proposition \ref{CGO:SolMagnetic} to choose solutions of the form 
   \begin{eqnarray*}
      & u = e^{-\rho/h} \big(\detg^{-1/4} c^{1/2} e^{i\Phi_1} a_0(x_1,r) b(\theta) + r_1\big), & \\
      & v = e^{\rho/h} \big(\detg^{-1/4} c^{1/2} e^{i\Phi_2} + r_2\big), & 
   \end{eqnarray*}
   where $\rho = x_1+ir$, $\dbar a_0 = 0$, and $\Phi = \Phi_1 + \Phi_2$ is a solution of the equation 
   \begin{equation} 
   \label{magnetic_phi_equation}
      \dbar \Phi + \frac{1}{2} (\tilde{A}_1 + i \tilde{A}_r) = 0 \quad \text{in $M$}.
   \end{equation}
   Here $\tilde{A} = A_1 - A_2$, and $\tilde{A}_1$ and $\tilde{A}_r$ are the components in the $x_1$ and $r$ coordinates. Inserting these solutions in   
   \eqref{magnetic_integral_identity}, multiplying both sides by $h$, and letting $h \to 0$ implies that 
   \begin{equation*}
      \lim_{h \to 0} \int_M \langle \tilde{A}, d\rho \rangle uv \,dV = 0.
   \end{equation*}
   Writing the integral in local coordinates gives 
   \begin{equation*}
      \int (\tilde{A}_1 + i \tilde{A}_r) e^{i\Phi} a_0(x_1,r) b(\theta) \,dx_1 \,dr \,d\theta = 0.
   \end{equation*}
   Varying $b$ leads to  
   \begin{equation*}
      \int_{\Omega_{\theta}} (\tilde{A}_1 + i \tilde{A}_r) e^{i\Phi} a_0(x_1,r) \,d\bar{\rho} \wedge d\rho = 0, \quad \text{for all $\theta$},
   \end{equation*}
   where $\Omega_{\theta} = \{(x_1,r) \in \R^2 \,;\, (x_1,r,\theta) \in M \}$ is identified with a domain in the complex plane, with complex variable $\rho$.       Integrating by parts and using \eqref{magnetic_phi_equation} gives 
   \begin{equation}  
   \label{magnetic_orthogonality_a0}
      \int_{\partial \Omega_{\theta}} e^{i\Phi} a_0 \,d\rho = 0.
   \end{equation}
   The arguments in \cite[Section 5]{DSFKSU}, see also \cite[Section 7]{KnSa}, then imply that $e^{i\Phi}|_{\partial \Omega_{\theta}} = F|_{\partial \Omega_{\theta}}$ where 
   $F \in C(\overline{\Omega_{\theta}})$ is a nonvanishing holomorphic function, and $F = e^G$ where $G \in C(\overline{\Omega_{\theta}})$ is holomorphic 
   and $G - i\Phi$ is constant on $\partial \Omega_{\theta}$. We choose 
   \begin{equation*}
      a_0 = G e^{-G} e^{i\lambda(x_1+ir)}
   \end{equation*}
   where $\lambda$ is a real number, and then \eqref{magnetic_phi_equation}, \eqref{magnetic_orthogonality_a0}, and integration by parts imply 
   \begin{equation*}
      \int_{\Omega_{\theta}} (\tilde{A}_1 + i \tilde{A}_r) e^{i\lambda(x_1+ir)} \,d\bar{\rho} \wedge d\rho = 0, \quad \text{for all $\theta$}.
   \end{equation*}
   We define  
   \begin{gather*}
      f(x') = \int e^{i\lambda x_1} \tilde{A}_1(x_1,x') \,dx_1, \\
      \alpha(x') = \sum_{j=2}^n \Big( \int e^{i\lambda x_1} \tilde{A}_j(x_1,x') \,dx_1 \Big) dx^j.
   \end{gather*}
   The identity above may be written as 
   \begin{equation*}
     \int e^{-\lambda r} \big[f(\gamma(r)) + i\alpha(\dot{\gamma}(r)) \big] \,dr = 0, \quad \text{for all $\theta$}.
   \end{equation*}
   The $r$-integrals are integrals over geodesics $\gamma$ in $D$. By varying the point $\omega$ in Proposition \ref{CGO:SolMagnetic} on $\partial D$ and using Theorem 
   \ref{thm:attenuated_injectivity}, we see for $\lambda$ small that $f = -\lambda p$ and $\alpha = -idp$ where $p \in C^{\infty}(D)$ and $p|_{\partial D} = 0$.
   The definition of $\alpha$ and analyticity of the Fourier transform imply that 
   \begin{equation*}
      \partial_k \tilde{A}_j - \partial_j \tilde{A}_k = 0, \quad j,k = 2,\ldots,n.
   \end{equation*}
   Also 
   \begin{equation*}
      \int e^{i\lambda x_1} (\partial_j \tilde{A}_1 - \partial_1 \tilde{A}_j)(x_1,x') \,dx_1 = \partial_j f + i\lambda \alpha_j = 0,
   \end{equation*}
   showing that $d\tilde{A} = 0$ in $M$ and the magnetic fields coincide.

   Since $M$ is convex and $\tilde{A}|_{\partial M} = 0$, one has $\tilde{A} = dp$ where $p \in C^{\infty}(M)$ can be chosen so 
   that $p|_{\partial M} = 0$. By a gauge transformation, we may assume that $A_1 = A_2$, and this $1$-form will be denoted by $A$. 
   The integral identity \eqref{magnetic_integral_identity} becomes 
   \begin{equation} 
   \label{magnetic_second_identity}
      \int_M (q_1-q_2) uv \,dV = 0,
   \end{equation}
   for any solutions of $\mathcal{L}_{g,A,q_1} u = 0$ and $\mathcal{L}_{g,-A,q_2} v = 0$. Use Proposition \ref{CGO:SolMagnetic} 
   and choose solutions 
   \begin{eqnarray*}
      & u = e^{-\rho/h} \big(\detg^{-1/4} c^{1/2} e^{i\Phi} e^{i\lambda(x_1+ir)} b(\theta) + r_1\big), & \\
      & v = e^{\rho/h} \big(\detg^{-1/4} c^{1/2} e^{-i\Phi} + r_2\big), & 
   \end{eqnarray*}
   where $\rho = x_1+ir$, and $\Phi$ satisfies 
   \begin{equation*}
      \dbar \Phi + \frac{1}{2} (A_1 + i A_r) = 0 \quad \text{in $M$}.
   \end{equation*}
   Letting $h \to 0$ in \eqref{magnetic_second_identity} gives 
   \begin{equation*}
      \int e^{i\lambda(x_1+ir)} (q_1-q_2) c(x_1,r,\theta) b(\theta) \,dx_1 \,dr \,d\theta = 0.
   \end{equation*}
   Proceeding as in the proof of Theorem \ref{Intro:MainThm} shows that $q_1 = q_2$.
\end{proof}
\end{subsection} 
\begin{subsection}{Recovering a conformal factor}

The results on the Schr\"odinger inverse problem can be used to recover a conformal factor from the DN map. 
Recall that we use the notation $\Lambda_g=\Lambda_{g,0}$ when the potential $q$ is zero.

\begin{proof}[Proof of Theorem \ref{Intro:ConfThm}]
It is enough to show that if $(M,g)$ is admissible and $c$ is smooth and positive, and if $\Lambda_{c g} = \Lambda_g$, then $c = 1$. We have $c|_{\partial M} = 1$ and $\partial_{\nu} c|_{\partial M} = 0$ by Proposition \ref{BDet:ThmLeU}, and then the assumption and Proposition \ref{Bdet:GaugeInv} imply 
\begin{equation*}
\Lambda_{c g,0} = \Lambda_{g,0} = \Lambda_{cg,q},
\end{equation*}
where $q = - \Delta_{g}(c^{\frac{n-2}{4}})/c^{\frac{n+2}{4}}$. We conclude from Theorem \ref{Intro:MainThm} that $q = 0$, so $\Delta_{g} (c^{\frac{n-2}{4}}) = 0$ in $M$. Since $c^{\frac{n-2}{4}} = 1$ on $\partial M$, uniqueness of solutions for the Dirichlet problem shows that $c \equiv 1$.
\end{proof}

\end{subsection}
\end{section}
%
%
\begin{section}{Attenuated ray transform}

For the uniqueness results in inverse problems, we have used that certain geodesic ray transforms are injective. If $(M,g)$ is a compact manifold with smooth boundary, geodesics can be parametrized by points on the unit sphere bundle $S M = \{ (x,\xi) \in TM \,;\, \abs{\xi} = 1 \}$. For $(x,\xi) \in SM$ let $\gamma_{x,\xi}(t)$ be the geodesic with $\gamma(0) = x$ and $\dot{\gamma}(0) = \xi$. We assume that $(M,g)$ is nontrapping, which means that the time $\tau(x,\xi)$ when $\gamma_{x,\xi}$ exits $M$ is always finite.

Given a smooth real function $a$ on $M$, the attenuated geodesic ray transform of a function $f$ is given by 
\begin{equation*}
I^a f(x,\xi) = \int_0^{\tau(x,\xi)} f(\gamma_{x,\xi}(t)) \exp\Big[ \int_0^t a(\gamma_{x,\xi}(s)) \,ds \Big] dt
\end{equation*}
for $(x,\xi) \in \partial_+ SM$. Here we use the sets of inward and outward pointing unit vectors 
\begin{equation*}
\partial_{\pm} S(M) = \{ (x,\xi) \in SM \,;\, x \in \partial M, \,\pm \langle \xi, \nu(x) \rangle < 0 \},
\end{equation*}
and $\nu$ is the outer unit normal vector to $\partial M$.

We will also need to integrate $1$-forms over geodesics. Let $f$ be a smooth function and $\alpha = \alpha_i \,dx^i$ a smooth $1$-form on $M$, and consider 
\begin{equation} \label{GRT:fform}
F(x,\xi) = f(x) + \alpha_i(x) \xi^i
\end{equation}
for $(x,\xi) \in SM$. The attenuated geodesic ray transform of $F$ is 
\begin{equation*}
I^a F(x,\xi) = \int_0^{\tau(x,\xi)} F(\gamma_{x,\xi}(t), \dot{\gamma}_{x,\xi}(t)) \exp\Big[ \int_0^t a(\gamma_{x,\xi}(s)) \,ds \Big] dt.
\end{equation*}
This transform always has a kernel: if $p$ is a smooth function on $M$ with $p|_{\partial M} = 0$, then a direct computation shows that 
\begin{equation*}
I^a(ap + dp(\xi)) = 0.
\end{equation*}
The main result in this section states that for simple manifolds and small attenuation, this is the only obstruction to injectivity.

\begin{thm} 
\label{thm:attenuated_injectivity}
Let $(M,g)$ be a compact simple manifold with smooth boundary. There exists $\varepsilon>0$ such that the following assertion holds for any smooth real function $a$ on $M$ satisfying $\abs{a} \leq \varepsilon$: If $F$ is given by \eqref{GRT:fform} and if 
\begin{equation*}
     I^a F(x,\xi) = 0
\end{equation*}
for all $(x,\xi) \in \partial_+ SM$, then $F = ap + dp(\xi)$ for some smooth function $p$ on $M$ with $p|_{\partial M} = 0$.
\end{thm}

Note that if $\alpha = 0$, this shows that any function $f$ whose attenuated ray transform vanishes must be identically zero. If $f = 0$ and $a \neq 0$ everywhere, then any $1$-form whose attenuated ray transform vanishes must be identically zero. Injectivity of the geodesic ray transform for functions and $1$-forms on simple manifolds in the case $a=0$ is well known \cite{An}, \cite{M}, \cite{Sh}. The injectivity for functions and small $a$ is proved in \cite{Sh}, \cite{Sh2} under conditions which involve a modified Jacobi equation or the size and curvature of $M$. We give a proof which works on simple manifolds.

We remark that the notation in this section is somewhat different from the other sections. For instance, we will denote by $\nabla$ the covariant derivative and more generally the horizontal derivative. The notation will be explained in more detail below.

\subsection{Preliminaries}

The proof will be based on energy estimates and a Pestov identity, which is the standard approach to such problems. First we need to recall the definition of horizontal (or semibasic) tensor fields on $TM$. These are $(p,q)$ tensor fields on $TM$ which have coordinate representations 
\begin{equation*}
u = (u_{i_1 \cdots i_p}^{j_1 \cdots j_q}) = u_{i_1 \cdots i_p}^{j_1 \cdots j_q} \frac{\partial}{\partial \xi^{j_1}} \otimes \cdots \otimes \frac{\partial}{\partial \xi^{j_q}} \otimes dx^{i_1} \otimes \cdots \otimes dx^{i_p}
\end{equation*}
with respect to coordinates $(x,\xi)$ on $TM$ associated to charts $x$ of $M$. The components transform in the same way as tensors on $M$ under changes of charts. Tensor fields on $M$ may be considered as $\xi$-constant horizontal tensor fields on $TM$. See \cite{Sh} for an invariant definition and other details on horizontal tensor fields.

For our purposes, it is sufficient to know that smooth functions on $TM$ are horizontal tensors of degree $(0,0)$, and that the horizontal and vertical derivatives defined by 
\begin{gather*}
(\overset{h}{\nabla} u)_{i_1 \cdots i_p i}^{j_1 \cdots j_q} = \overset{h}{\nabla}_i u_{i_1 \cdots i_p}^{j_1 \cdots j_q} := \tilde{\nabla}_i u_{i_1 \cdots i_p}^{j_1 \cdots j_q} - \Gamma_{ik}^l \xi^k \frac{\partial}{\partial \xi^l} u_{i_1 \cdots i_p}^{j_1 \cdots j_q}, \\
(\overset{v}{\nabla} u)_{i_1 \cdots i_p i}^{j_1 \cdots j_q} = \overset{v}{\nabla}_i u_{i_1 \cdots i_p}^{j_1 \cdots j_q} := \frac{\partial}{\partial \xi^i} u_{i_1 \cdots i_p}^{j_1 \cdots j_q},
\end{gather*}
are invariantly defined operators which map horizontal $(p,q)$ tensors to horizontal $(p+1,q)$ tensors. Here 
\begin{equation*}
\tilde{\nabla}_i u_{i_1 \cdots i_p}^{j_1 \cdots j_q} = \frac{\partial}{\partial x^i} u_{i_1 \cdots i_p}^{j_1 \cdots j_q} + \sum_{r=1}^q u_{i_1 \cdots i_p}^{j_1 \cdots s \cdots j_q} \Gamma_{is}^{j_r} - \sum_{r=1}^p u_{i_1 \cdots s \cdots i_p}^{j_1 \cdots j_q} \Gamma_{i i_r}^s.
\end{equation*}
Thus $\tilde{\nabla}$ acts in the same way as the usual covariant derivative.

Below, we will work with smooth functions and tensors on $SM$. Let $p: TM \msetminus \{0\} \to SM$ be the map $(x,\xi) \mapsto (x,\xi/\abs{\xi})$. The horizontal and vertical derivatives on $SM$ are defined by 
\begin{equation*}
\nabla u = \overset{h}{\nabla}(u \circ p)|_{SM}, \quad \partial u = \overset{v}{\nabla}(u \circ p)|_{SM}.
\end{equation*}
Also $\nabla_i u = \overset{h}{\nabla}_i(u \circ p)|_{SM}$ and $\partial_i u = \overset{v}{\nabla}_i(u \circ p)|_{SM}$. If $u$ is horizontal tensor field on $SM$, then $\nabla u$ and $\partial u$ are also horizontal tensor fields. For a smooth function on $SM$, these derivatives are given by 
\begin{align*}
\nabla_i u &= \frac{\partial}{\partial x^i} (u(x,\xi/\abs{\xi})) - \Gamma_{ik}^l \xi^k \partial_l u, \\
\partial_i u &= \frac{\partial}{\partial \xi^i} (u(x,\xi/\abs{\xi})).
\end{align*}
We mention the following formulas
\begin{gather*} 
\nabla g=0,\quad \nabla \xi =0,\quad \partial _{j}\xi ^{i}= \delta
_{j}^{i} - \xi^i \xi_j, \\
\lbrack \nabla_i ,\partial_j \rbrack =0,\quad \lbrack \partial
_{i},\partial _{j}\rbrack = \xi_i \partial_j - \xi_j \partial_i, \\
\lbrack \nabla _{i},\nabla _{j}\rbrack u = -R_{ijkl} \xi^k \partial^l u,
\end{gather*}
where $R$ is the curvature tensor and $u$ is a scalar function. We write 
\begin{equation*}
\partial^i u = g^{ij} \partial_j u, \quad \langle \partial u, \partial v \rangle = \partial^i u \partial_i v, \quad \abs{\partial u}^2 = \partial^i u \partial_i u.
\end{equation*}

Let $H$ be the geodesic vector field on $SM$ that generates geodesic flow. In local 
coordinates 
\[
Hu(x,\xi)=\xi^i \frac{\partial}{\partial x^i} (u(x,\xi/\abs{\xi}))-\Gamma
^l_{ik}(x)\xi^i \xi^k \partial_l u, \quad \xi\in S_x,
\]
where $u$ is a smooth function on $SM$. We may apply the operator $H$ to horizontal tensor fields on $SM$ by
defining $Hu = \xi^i \nabla_i u$.

\subsection{Ray transform of functions}
We now consider the boundary value problem for transport equation:
\[
Hu+au=-f, \quad u|_{\partial_-S(M)}=0,
\]
where the absorption $a$ and $f$ are smooth functions on $(M,g)$. 
The solution is given by
\begin{equation*}
u^f_a(x,\xi)=\int_0^{\tau
(x,\xi)}f(\gamma_{x,\xi}(t))\exp\Big[\int_0^t a(\gamma_{x,\xi}(s)) \,ds\Big] dt,
\end{equation*}
The trace
\[
I^af=u^f_a|_{\partial _{+}S(M)}
\]
is the attenuated geodesic ray transform of the function $f$. It is natural to define $\tau|_{\partial_{-} SM} = 0$, and then $u^f_a$ indeed vanishes on $\partial_{-} SM$.

We will prove that $f$ is uniquely determined by $I^a f$ under the following assumption. Let $\gamma=\gamma_{x,\xi}$, $(x,\xi)\in \partial_+ SM$, be an arbitrary geodesic, and consider the quadratic form 
\begin{equation} \label{modified_jacobi}
E_\gamma^a(X)=\int_0^{\tau} \big(|DX|^2-\langle R_{\gamma}X,X \rangle-a^2|X|^2\big)(t)\,dt
\end{equation}
where $\tau = \tau(x,\xi)$, $X(t)$ is a vector field on $\gamma$ belonging to the space 
\begin{equation*}
H^1_0(\gamma) = \{X\in H^1([0,\tau];T(\gamma)) \,;\, X(0)=X(\tau)=0\},
\end{equation*}
$D$ is the covariant derivative along $\gamma$, and $R_{\gamma} X = R(X,\dot{\gamma}) \dot{\gamma}$. We assume for all geodesics $\gamma$ the positive definiteness of this quadratic form, 
\begin{align} \label{modified_index}
   \left\{ 
   \begin{gathered}
      E_\gamma^a(X) \geq 0 \quad \text{whenever} \quad X \in H^1_0(\gamma), \\
      E_\gamma^a(X) = 0 \quad \text{iff} \quad X = 0.
   \end{gathered}
   \right.   
\end{align}
If \eqref{modified_index} holds we say that any geodesic has no conjugate points with respect to \eqref{modified_jacobi}. If $a = 0$ we obtain the usual index form $E_{\gamma} = E_{\gamma}^0$. Then clearly there are no conjugate points in the usual sense if there are none with respect to \eqref{modified_jacobi}.

\begin{prop} \label{prop:attenuated_function}
Let $(M,g)$ be compact and $\partial M$ strictly convex. Assume that any
geodesic has no conjugate points with respect to \eqref{modified_jacobi}.
Then any smooth function $f$ on the manifold $(M,g)$ is uniquely
determined by its attenuated geodesic X-ray transform.
\end{prop}
\begin{proof}
Let $I^af=0$. We will assume that $u=u^f_a$ is smooth on $SM$ (otherwise one can work in a slightly smaller manifold than $M$, and pass to the limit using the smoothness properties of $\tau$ as in \cite{Sh}). The function $u$ satisfies
$$\partial Hu+a\partial u=0$$ and therefore 
\begin{equation} \label{dhu_equation}
|\partial Hu|^2 = a^2|\partial u|^2.
\end{equation}
Using the formulas for $\nabla$ and $\partial$, it is not difficult to check the following
identity (valid for any $u \in C^{\infty}(SM)$):
\begin{equation} \label{pestov_identity}
|\partial Hu|^2=|H\partial u|^2 + \delta V + \theta W - R(\partial u,\xi,\xi,\partial u),
\end{equation}
where $\delta$ and $\theta$ are the vertical and horizontal divergences,
\[
\delta X=\nabla_iX^i, \quad \theta X=\partial_iX^i,
\]
and $V$ and $W$ are defined by 
\[
V^i = \langle \partial u, \nabla u \rangle \xi^i - (Hu) \partial^i u, \quad W^i = (Hu) \nabla^i u.
\]

From \eqref{dhu_equation} and \eqref{pestov_identity} we obtain 
$$
|H\partial u|^2-R(\partial u,\xi,\xi,\partial u) -a^2|\partial
u|^2 + \delta V + \theta W =0.
$$
Now integrate this equality over the manifold $SM$. Before
this we recall the integration formulas (see \cite{Sh}):
\begin{gather*}
\int_{SM} v \,d(SM) = \int_{M} \,dM \int_{S_x} v \,dS_x, \\
\int_{S_x}\theta X \,dS_x=(n-1)\int_{S_x} \langle X,\xi \rangle \,dS_x, \\
\nabla\int_{S_x}u \,dS_x=\int_{S_x}\nabla u \,dS_x,
\end{gather*}
where $v$ is a scalar, $X$ is a horizontal vector field and $u$ is horizontal tensor field,
\begin{equation*}
\int_{M} \delta X \,dM = \int_{\partial M} \langle X, \nu \rangle \,d(\partial M),
\end{equation*}
where $X$ is a vector field on $(M,g)$. In these formulas the volume
forms of corresponding manifolds are naturally defined using  the
metric $g$. After integration we have
\begin{multline} \label{pestov_integrated}
\int_{SM} \big(|H\partial u|^2-R(\partial u,\xi,\xi,\partial u)-a^2|\partial u|^2\big) \,d(SM) \\
 +(n-1)\int_{SM}|Hu|^2 \,d(SM) =0.
\end{multline}
We used the fact that $\langle V, \nu \rangle$ vanishes on $\partial (SM)$ since $u|_{\partial(SM)} = 0$.

We next show that \eqref{modified_index} implies
\[
E^a(Y)=\int_{SM}\big(|HY|^2-R(Y,\xi,\xi,Y)-a^2|Y|^2\big) \,d(SM)\geq 0,
\]
for any horizontal vector field $Y\in C^{\infty}(SM;TM)$ with $Y|_{\partial (SM)}=0$, and that the equality holds iff
$Y=0$. Santal\'o's formula (see \cite{Sh}) states
\[
\int_{SM} v \,d(SM)=
-\int_{\partial_+S(M)} \int_0^{\tau(x,\xi)} v(\gamma_{x,\xi}(t), \dot{\gamma}_{x,\xi}(t)) \langle \xi, \nu \rangle \,dt \,d(\partial (SM)),
\]
for $v \in C^{\infty}(SM)$. Let 
\begin{equation*}
X(x,\xi,t)=Y(\gamma_{x,\xi}(t),\dot{\gamma}_{x,\xi}(t)),
\end{equation*}
which implies $HY(\gamma_{x,\xi}(t),\dot{\gamma}_{x,\xi}(t)) = DX(x,\xi,t)$. Then we have
$$E^a(Y)=-\int_{\partial_+S(M)} E^a_{\gamma_{x,\xi}}(X)(x,\xi) \langle \xi,\nu(x) \rangle \,d(\partial (SM)) \geq 0.$$
Equality holds iff $Y = 0$. We have from \eqref{pestov_integrated} that $Hu=0$, which implies $u=0$ and
$f=0$.
\end{proof}

\subsection{Ray transform of $1$-forms}

Let $f$ be a smooth function and $\alpha = \alpha_i(x) dx^i$ a smooth $1$-form in $M$, and consider the boundary value problem 
\begin{equation*}
Hu + au = -F, \quad u|_{\partial_- SM} = 0,
\end{equation*}
where $a$ is a smooth function on $M$ and $F$ is as in \eqref{GRT:fform}. The solution $u = u_a^F$ is given by 
\begin{equation*}
u^F_a(x,\xi)=\int_0^{\tau(x,\xi)} F(\gamma_{x,\xi}(t),\dot{\gamma}_{x,\xi}(t))\exp\Big[\int_0^t a(\gamma_{x,\xi}(s)) \,ds\Big] \,dt,
\end{equation*}
and the trace 
\begin{equation*}
I^a F = u_a^F|_{\partial_+ SM}
\end{equation*}
is the attenuated geodesic X-ray transform of $F$.

\begin{prop} \label{prop:attenuated_oneform}
Let $(M,g)$ be compact and $\partial M$ strictly convex, and suppose that any geodesic has no conjugate points with respect to \eqref{modified_jacobi}. If $I^a F = 0$, then $F = ap + dp(\xi)$ for some smooth function $p$ on $M$ which vanishes on $\partial M$.
\end{prop}
\begin{proof}
We follow the proof of Proposition \ref{prop:attenuated_function}. Again assume that $u$ is smooth in $SM$. Then $u$ satisfies 
\begin{equation*}
\partial Hu + a \partial u = -\partial F,
\end{equation*}
and 
\begin{equation*}
\abs{\partial Hu}^2 = a^2 \abs{\partial u}^2 + 2a \langle \partial u, \partial F \rangle + \abs{\partial F}^2.
\end{equation*}
The identity \eqref{pestov_identity} then implies 
\begin{equation*}
\abs{H\partial u}^2 - R(\partial u,\xi,\xi,\partial u) - a^2\abs{\partial u}^2 + \delta V + \theta W - 2a \langle \partial u, \partial F \rangle - \abs{\partial F}^2 = 0.
\end{equation*}
Integrating this over $SM$, it follows that 
\begin{equation} \label{pestov_integrated_oneform}
E^a(\partial u) + \int_{SM} \big[(n-1)\abs{Hu}^2 - 2a \langle \partial u, \partial F \rangle - \abs{\partial F}^2\big] \,d(SM) = 0.
\end{equation}
Using the integration formula 
\begin{equation*}
\int_{S_x} \partial_i v \,dS_x = (n-1) \int_{S_x} v \xi_i \,dS_x
\end{equation*}
and the identities 
\begin{gather*}
\abs{Hu}^2 = (au)^2 + 2auF + F^2, \\
\partial_i \partial^i F = -(n-1)\alpha_i \xi^i,
\end{gather*}
we easily obtain from \eqref{pestov_integrated_oneform} that 
\begin{equation*}
E^a(\partial u) + (n-1) \int_{SM} (au + f)^2 \,d(SM) = 0.
\end{equation*}
The assumption \eqref{modified_index} and Santal\'o's formula imply $E^a(\partial u) \geq 0$, so from the last equality we obtain $E^a(\partial u) = 0$ and $au+f = 0$. Now $E^a(\partial u) = 0$ implies $\partial u = 0$, so $u = u(x)$ and 
\begin{equation*}
\langle du, \xi \rangle + au = -f -\langle \alpha, \xi \rangle, \quad x \in M, \,\xi \in S_x.
\end{equation*}
The claim follows.
\end{proof}

It remains to show that simple manifolds satisfy the condition of Proposition \ref{prop:attenuated_oneform}.

\begin{proof}[Proof of Theorem \ref{thm:attenuated_injectivity}]
Let $(M,g)$ be a compact simple manifold with smooth boundary. We need to show that there is $\varepsilon > 0$ such that 
\begin{equation} \label{index_strictly_positive}
E_{\gamma_{x,\xi}}(X) \geq \varepsilon \int_0^{\tau(x,\xi)} \abs{X}^2 \,dt,
\end{equation}
for all $(x,\xi) \in \overline{\partial_+ SM}$ and for all $X \in H^1_0(\gamma_{x,\xi})$. If \eqref{index_strictly_positive} holds and $\abs{a} < \sqrt{\varepsilon}$ then any geodesic on $(M,g)$ has no conjugate points with respect to \eqref{modified_jacobi}, and Proposition \ref{prop:attenuated_oneform} implies the desired result.

Let first $(x,\xi) \in \partial_+ SM$, and consider the unbounded operator on $L^2(\gamma_{x,\xi})$ with domain $H^2 \cap H^1_0(\gamma_{x,\xi})$, given by 
\begin{equation*}
L_{\gamma_{x,\xi}}: X \mapsto -D^2 X - R(X,\dot{\gamma}_{x,\xi}) \dot{\gamma}_{x,\xi}.
\end{equation*}
This operator is self-adjoint and has discrete spectrum, which lies in $(0,\infty)$ since the corresponding quadratic form $E_{\gamma_{x,\xi}}$ is positive definite. Therefore 
\begin{equation*}
E_{\gamma_{x,\xi}}(X) \geq \lambda_1(x,\xi)  \int_0^{\tau(x,\xi)} \abs{X}^2 \,dt, \quad X \in H^1_0(\gamma_{x,\xi}),
\end{equation*}
where $\lambda_1(x,\xi) > 0$ is the smallest eigenvalue. The coefficients of $L$ depend smoothly on $(x,\xi)$ and $\tau$ is smooth and positive in $\partial_+ SM$, so it is not hard to see that \eqref{index_strictly_positive} holds in a neighborhood of any fixed point in $\partial_+ SM$.

For tangential directions we use the Poincar\'e inequality 
\begin{equation*}
\int_0^L \abs{u(t)}^2 \,dt \leq \frac{L^2}{\pi^2} \int_0^L \abs{\dot{u}(t)}^2 \,dt, \quad u \in H^1_0([0,L]),
\end{equation*}
where the constant $L^2/\pi^2$ is optimal \cite{PW}. If $(x,\xi) \in S(\partial M)$ and $\delta > 0$, then by continuity of $\tau$ there is a neighborhood $U$ of $(x,\xi)$ in $\overline{\partial_+ SM}$ such that $\tau(y,\eta) \leq \delta$ in that neighborhood. Choosing $\delta$ small enough, the Poincar\'e inequality implies 
\begin{multline*}
E_{\gamma_{y,\eta}}(X) = \int_0^{\tau(y,\eta)} (\abs{DX}^2 - \langle R_{\gamma_{y,\eta}} X, X \rangle) \,dt \\
 \geq \int_0^{\tau(y,\eta)} \Big( \frac{\pi^2}{\tau(y,\eta)^2} \abs{X}^2 - \langle R_{\gamma_{y,\eta}} X, X \rangle \Big) \,dt \geq \frac{\pi^2}{2\delta^2} \int_0^{\tau(y,\eta)} \abs{X}^2 \,dt
\end{multline*}
whenever $(y,\eta) \in U$ and $X \in H^1_0(\gamma_{y,\eta})$. This shows \eqref{index_strictly_positive} near any point of $S(\partial M)$. It follows that for some $\varepsilon > 0$, \eqref{index_strictly_positive} holds on the compact set $\overline{\partial_+ SM}$.
\end{proof}

\end{section}
%
%
\begin{section}{Boundary determination}

To deal with the inverse problems we are interested in, we need to use the fact that the
DN map determines the Taylor expansions at the boundary of the different quantities involved. In the case
of the Laplace-Beltrami operator, the relevant result is proved in \cite{LeU} and is as follows.

\begin{prop}
\label{BDet:ThmLeU}
   Let $(M,g_1)$ and $(M,g_2)$ be compact manifolds with smooth boundary, with dimension $n \geq 3$. If $\Lambda_{g_1}=\Lambda_{g_2}$,
   then the Taylor series of $g_1$ and $g_2$ in boundary normal coordinates are equal at each point on the boundary.
\end{prop}

In this section we extend the previous result to the case where electric and magnetic potentials are present.
First we need to consider the gauge invariance of the DN map.

\begin{prop}
\label{Bdet:GaugeInv}
   Let $(M,g)$ be a compact manifold with boundary, and let $A$ be a smooth 1-form and $q$ a 
   smooth function on $M$. If $c$ and $\psi$ are smooth functions such that 
       $$ c>0, \quad c|_{\d M}=1, \quad \d_{\nu}c|_{\d M}=0, \quad \psi|_{\d M}=0, $$
   then we have 
   \begin{align*}
      \Lambda_{g,A,q}=\Lambda_{c^{-1}g,A+d\psi,c(q-q_c)}
   \end{align*}
   where $q_c=c^{\frac{n-2}{4}}\Delta_g \big(c^{-\frac{n-2}{4}}\big)$.
\end{prop}
\begin{proof}
   It follows from a direct computation that 
   \begin{align*}
   		c^{\frac{n+2}{4}} \mathcal{L}_{g,A,q} (c^{-\frac{n-2}{4}} u) &= \mathcal{L}_{c^{-1} g,A,c(q-q_c)} u, \\
   		e^{-i\psi} \mathcal{L}_{g,A,q} (e^{i\psi} u) &= \mathcal{L}_{g,A+d\psi,q} u.
   \end{align*}
   Let $f \in C^{\infty}(\partial M)$, and let $u$ be the solution of $\mathcal{L}_{g,A,q} u = 0$ in $M$ which satisfies $u|_{\partial M} = f$. 
   If 
   \begin{equation*}
   		(\tilde{g},\tilde{A},\tilde{q}) = (c^{-1} g,A+d\psi,c(q-q_c))
   \end{equation*}
   and if $\tilde{u} = c^{\frac{n-2}{4}} e^{-i\psi} u$, we have 
   $\mathcal{L}_{\tilde{g},\tilde{A},\tilde{q}} \tilde{u} = 0$ in $M$ and $\tilde{u}|_{\partial M} = f$. Then 
   \begin{equation*}
   		\Lambda_{\tilde{g},\tilde{A},\tilde{q}} f = d_{\tilde{A}} \tilde{u}(\nu_{\tilde{g}}) |_{\partial M} = d_{A+d\psi}(e^{-i\psi} u)(\nu_g)|_{\partial M} = \Lambda_{g,A,q} f,
   \end{equation*}
   by using the boundary values of $c$ and $\psi$ and the fact that $\nu_{\tilde{g}} = \nu_g$.
\end{proof}   

\begin{rem}
   The conformal gauge invariance is related to the fact that the conformal Laplace-Beltrami operator $\tilde{\Delta}_g$ on $(M,g)$, defined by 
       $$ \tilde{\Delta}_g = \Delta_g - \frac{n-2}{4(n-1)} \mathop{\rm Scal}\nolimits_g, $$
   transforms under a conformal change of metrics by 
       $$ \tilde{\Delta}_{c g}u = c^{-\frac{n+2}{4}} \tilde{\Delta}_g\big(c^{\frac{n-2}{4}} u\big). $$
\end{rem}

We use the notation $f_1 \simeq f_2$ to denote that $f_1$ and $f_2$ have the same Taylor series. Our main boundary determination result is as follows.

\begin{thm}
\label{Bdet:ThmBDet}
   Let $(M,g_1)$ and $(M,g_2)$ be compact manifolds with boundary, of dimension $n \geq 3$, 
   and let $A_1,A_2$ be two smooth $1$-forms and $q_1,q_2$ two smooth functions in $M$. 
   If $\Lambda_{g_1,A_1,q_1} = \Lambda_{g_2,A_2,q_2}$ and if $p \in \partial M$, then there exist smooth positive functions 
   $c_j$ with 
      $$ c_j|_{\partial M} = 1 \quad \textrm{and} \quad \partial_{\nu} c_j|_{\partial M} = 0, $$
   and smooth functions $\psi_j$ with $\psi_j|_{\partial M} = 0$, such that the gauge transformed coefficients 
   \begin{equation*}
      (\tilde{g}_j,\tilde{A}_j,\tilde{q}_j) = \big(c_j^{-1} g_j, A_j + d\psi_j, c_j(q_j - d_j^{-1} \Delta_{g_j} d_j)\big)
   \end{equation*}
   with $d_j = c_j^{-\frac{n-2}{4}}$, satisfy in boundary normal coordinates at $p$ 
   \begin{equation*}
      \tilde{g}_1 \simeq \tilde{g}_2, \quad \tilde{A}_1 \simeq \tilde{A}_2, \quad \tilde{q}_1 \simeq \tilde{q}_2.
   \end{equation*}
    Furthermore, if $g_1 \simeq g_2$ in boundary normal coordinates on all of $\partial M$, then $\tilde{A}_1 \simeq \tilde{A}_2$ and 
    $q_1 \simeq q_2$ on $\partial M$.
\end{thm}

Boundary normal coordinates $(x',x_n)$ at a boundary point $p$ are such that $x'$ parametrizes $\partial M$ near $p$ and $x_n$ is the distance to the boundary along unit speed geodesics normal to $\partial M$. See \cite{LeU} for more details.

We prove Theorem \ref{Bdet:ThmBDet} by showing that $\Lambda_{g,A,q}$ is a pseudodifferential operator whose symbol determines the 
boundary values of the coefficients. This method was used in \cite{SU2} for the conductivity equation, in \cite{LeU} for the Laplace-Beltrami operator, and in \cite{NSU} for the magnetic Schr\"odinger operator with Euclidean metric. We follow the proof in \cite{LeU}.

The first issue to consider is the gauge invariance in the operator $\Lambda_{g,A,q}$. This will be dealt with by normalizing the coefficients $(g,A,q)$ in a way which is suitable for boundary determination. Fix a point $p \in \partial M$ and boundary normal coordinates $(x',x_n)$ near $p$. In these coordinates $\partial M$ corresponds to $\{x_n = 0\}$ and 
  $$ g = g_{\alpha \beta} \,dx^{\alpha} \otimes dx^{\beta} + dx^n \otimes dx^n. $$ 
We use the convention that Greek indices run from $1$ to $n-1$ and Roman indices from $1$ to $n$.

\begin{lem} 
\label{lemma:boundary_normalization}
Let $(g,A,q)$ be smooth coefficients in $M$ and $p \in \partial M$. There exist a positive smooth function $c$ with 
$c|_{\partial M} = 1$ and $\partial_{\nu} c|_{\partial M} = 0$, and a smooth function $\psi$ with $\psi|_{\partial M} = 0$, 
such that in boundary normal coordinates near $p$ the quantities $\tilde{g} = c^{-1} g$ and $\tilde{A} = A + d\psi$ satisfy 
\begin{align}
   \partial_n^j \tilde{A}_n(x',0) &= 0, \quad j \geq 0,  \label{a_normalization} \\
   \partial_n^j (\tilde{g}_{\alpha \beta} \partial_n \tilde{g}^{\alpha \beta})(x',0) &= 0, \quad j \geq 1.  \label{g_normalization}
\end{align}
\end{lem}
\begin{proof}
One can find a smooth function $\psi$ near $p$ with $\psi(x',0) = 0$ and $\partial_n^{j+1} \psi(x',0) = -\partial_n^j A_n(x',0)$ for $j \geq 0$. Extending this in a suitable way, one obtains $\psi \in C^{\infty}(M)$ with $\psi|_{\partial M} = 0$ such that $\tilde{A} = A+d\psi$ will satisfy \eqref{a_normalization}.

Further, we construct a smooth function $c$ near $p$ which satisfies $c(x',0) = 1$, $\partial_n c(x',0) = 0$, and 
\begin{equation*}
\partial_n^j (\log\,\det(cg^{\alpha \beta}))(x',0) = 0, \quad j \geq 2.
\end{equation*}
In fact, one may take $c = e^{\mu}$ where $\mu(x',0) = \partial_n \mu(x',0) = 0$ and $\partial_n^j \mu(x',0) = -\frac{1}{n-1} \partial_n^j (\log\,\det(g^{\alpha \beta}))(x',0)$ for $j \geq 2$. There is an extension of $c$ to a positive function $c \in C^{\infty}(M)$ with $c|_{\partial M} = 1$ and $\partial_{\nu} c|_{\partial M} = 0$. Since 
\begin{equation*}
\partial_n (\log\,\det(\tilde{g}^{\alpha \beta})) = \tilde{g}_{\alpha \beta} \partial_n \tilde{g}^{\alpha \beta},
\end{equation*}
one also has the second condition \eqref{g_normalization}.
\end{proof}

Replacing coefficients $(g,A,q)$ by the gauge equivalent coefficients $(c^{-1} g, A+d\psi,c(q-d^{-1} \Delta_g d))$ does not affect the DN map. Therefore, below we will assume that \eqref{a_normalization} and \eqref{g_normalization} are valid. Note that boundary normal coordinates for $g$ are also boundary normal coordinates for any conformal multiple of $g$, if the conformal factor is $1$ on $\partial M$.

The next aim is to write $\Lambda_{g,A,q}$ as a pseudodifferential operator and to compute the symbol in a small neighborhood of $p$. Here we use the usual (not semiclassical) pseudodifferential calculus and left quantization, so that a symbol $r(x,\xi)$ in $T^* \R^n$ corresponds to the operator 
\begin{equation*}
R f(x) = (2\pi)^{-n} \int_{\R^n} \int_{\R^n} e^{i(x-y) \cdot \xi} r(x,\xi) f(y) \,dy \,d\xi.
\end{equation*}
We denote by $p \sim \sum p_j$ the asymptotic sum of symbols, see \cite{H3} for these basic facts.

\begin{lem} 
\label{lemma:dnmap_psdo}
$\Lambda_{g,A,q}$ is a pseudodifferential operator of order $1$ on $\partial M$. Its full symbol (in left quantization) in boundary normal coordinates near $p$ is $-b \sim -\sum_{j \leq 1} b_j$, where $b_j$ are given in \eqref{b1_formula} -- \eqref{bmm1_formula}.
\end{lem}
\begin{proof}
In the $x$ coordinates, one has 
\begin{equation*}
\mathcal{L}_{g,A,q} = -\Delta_g + 2 g^{jk} A_j D_k + G,
\end{equation*}
where $G = \detg^{-1/2} D_j(\detg^{1/2} g^{jk} A_k) + g^{jk} A_j A_k + q$ and 
 $$\detg = \det(g_{jk}) = \det(g_{\alpha \beta}).$$
From \eqref{a_normalization} we know that $\partial^K A_n(x',0) = 0$ for all multi-indices $K \in \N^n$.

One would like to have $A_n = 0$ also inside $M$. To achieve this, we introduce as in \cite{NSU} the conjugated operator 
\begin{equation*}
\mathcal{M} = e^{-ih} \mathcal{L}_{g,A,q} e^{ih},
\end{equation*}
where $h(x) = - \int_0^{x_n} A_n(x',s) \,ds$. Note that $\partial^K h(x',0) = 0$. Writing $\tilde{A}_j = A_j + \partial_j h$, we obtain $\tilde{A}_n = 0$ and 
\begin{equation*}
\mathcal{M} = -\Delta_g + 2 g^{\alpha \beta} \tilde{A}_{\alpha} D_{\beta} + \tilde{G},
\end{equation*}
where $\tilde{G} = \detg^{-1/2} D_{\alpha}(\detg^{1/2} g^{\alpha \beta} \tilde{A}_{\beta}) + g^{\alpha \beta} \tilde{A}_{\alpha} \tilde{A}_{\beta} + q$. We then have 
\begin{equation*}
\mathcal{M} = D_n^2 + i E(x) D_n + Q_2(x,D_{x'}) + Q_1(x,D_{x'}) + 2g^{\alpha \beta} \tilde{A}_{\alpha} D_{\beta} + \tilde{G},
\end{equation*}
with $E$, $Q_1$, and $Q_2$ given by 
\begin{gather*}
E(x) = \frac{1}{2} g_{\alpha \beta} \partial_n g^{\alpha \beta}, \\
Q_2(x,D_{x'}) = g^{\alpha \beta} D_{\alpha} D_{\beta}, \\
Q_1(x,D_{x'}) = -i(\frac{1}{2} g^{\alpha \beta} \partial_{\alpha}(\log\,\detg) + \partial_{\alpha} g^{\alpha \beta}) D_{\beta}.
\end{gather*}

As in \cite{LeU}, one would like to have a factorization 
\begin{multline} \label{m_factorization}
\mathcal{M} = (D_n + iE(x) - iB(x,D_{x'}))(D_n + iB(x,D_{x'})) \\
 \text{modulo a smoothing operator},
\end{multline}
where $B$ is a pseudodifferential operator of order $1$ with symbol $b(x,\xi')$. Using left quantization for symbols, \eqref{m_factorization} implies on the level of symbols that 
\begin{equation} \label{b_symbol_equation}
\partial_n b - E(x) b + \sum_K \frac{\partial_{\xi'}^K b D_{x'}^K b}{K!} = q_2 + q_1 + 2g^{\alpha \beta} \tilde{A}_{\alpha} \xi_{\beta} + \tilde{G}  \mod S^{-\infty}.
\end{equation}
Let $b \sim \sum_{j \leq 1} b_j$ where $b_j(x,\xi')$ is homogeneous of order $j$ in $\xi'$. Inserting this in \eqref{b_symbol_equation} and collecting terms with the same order of homogeneity, one obtains $b_j$ as follows:
\begin{align}
b_1 &= -\sqrt{q_2}, \label{b1_formula} \\
b_0 &= \frac{1}{2b_1}(-\partial_n b_1 + E b_1 - \nabla_{\xi'} b_1 \cdot D_{x'} b_1 + q_1 + 2g^{\alpha \beta} \tilde{A}_{\alpha} \xi_{\beta}), \label{b0_formula} \\
b_{-1} &= \frac{1}{2b_1}(-\partial_n b_0 + E b_0 - \sum_{\substack{0 \leq j,k \leq 1 \\ j+k=|K|}} \frac{\partial_{\xi'}^K b_j D_{x'}^K b_k}{K!} + \tilde{G}), \label{bm1_formula} \\
b_{m-1} &= \frac{1}{2b_1}(-\partial_n b_m + E b_m - \sum_{\substack{m \leq j,k \leq 1 \\j+k-|K|=m}} \frac{\partial_{\xi'}^K b_j D_{x'}^K b_k}{K!}) \ \  (m \leq -1). \label{bmm1_formula} 
\end{align}
With these choices, $b \in S^1$ and \eqref{m_factorization} is valid. By the argument in \cite[Proposition 1.2]{LeU}, one has 
\begin{equation*}
\Lambda_{g,A,q} f(x') = -B(x',0,D_{x'}) f(x') + Rf(x')
\end{equation*}
where $R$ is a smoothing operator.
\end{proof}

From the symbol of the DN map, one can recover the following information on the coefficients.

\begin{lem} 
\label{lemma:boundary_quantities}
The knowledge of $\Lambda_{g,A,q}$ determines on $\{x_n = 0\}$ the quantities 
\begin{equation} \label{boundary_quantities}
g^{\alpha \beta}, \,\partial_n g^{\alpha \beta}, \,\partial^K A_{\alpha}, \,\partial^K l^{\alpha \beta}.
\end{equation}
Here $K \in \N^n$ is any multi-index and 
\begin{equation*}
l^{\alpha \beta} = \frac{1}{4} \partial_n k^{\alpha \beta} + q g^{\alpha \beta},
\end{equation*}
with $k^{\alpha \beta} = \partial_n g^{\alpha \beta} - (g_{\gamma \delta} \partial_n g^{\gamma \delta}) g^{\alpha \beta}$.
\end{lem}
\begin{proof}
By Lemma \ref{lemma:dnmap_psdo}, $\Lambda_{g,A,q}$ determines $b|_{x_n=0}$ and each $b_j|_{x_n = 0}$. From $b_1|_{x_n=0}$ one recovers $g^{\alpha \beta}|_{x_n=0}$. If $T(g^{\alpha \beta})$ denotes any linear combination of tangential derivatives of $g^{\alpha \beta}$, one has for $x_n = 0$ 
\begin{align*}
b_0 &= \frac{1}{2b_1} (-\partial_n b_1 + E b_1 + 2 g^{\alpha \beta} \tilde{A}_{\alpha} \xi_{\beta}) + T(g^{\alpha \beta}) \\
 &= - \frac{\partial_n q_2}{4 q_2} + \frac{1}{2} E - \frac{1}{\sqrt{q_2}} g^{\alpha \beta} \tilde{A}_{\alpha} \xi_{\beta} + T(g^{\alpha \beta}) \\
 &= -\frac{1}{4} k^{\alpha \beta} \omega_{\alpha} \omega_{\beta} - g^{\alpha \beta} \tilde{A}_{\alpha} \omega_{\beta} + T(g^{\alpha \beta}),
\end{align*}
where $\omega = \xi'/|\xi'|_g$. Evaluating at $\pm \omega$ and varying $\omega$ one recovers $\tilde{A}_{\alpha}$ and $k^{\alpha \beta}$, and consequently also $\partial_n g^{\alpha \beta}$, for $x_n = 0$.

Moving to $b_{-1}$, one has for $x_n = 0$ 
\begin{align*}
b_{-1} &= \frac{1}{2b_1} (-\partial_n b_0 + E b_0 + \tilde{G}) + T(g^{\alpha \beta}, \partial_n g^{\alpha \beta}, \tilde{A}_{\alpha}) \\
 &= \frac{1}{2b_1} ( l^{\alpha \beta} \omega_{\alpha} \omega_{\beta} + g^{\alpha \beta}(\partial_n \tilde{A}_{\alpha}) \omega_{\beta}) + T(g^{\alpha \beta}, \partial_n g^{\alpha \beta}, \tilde{A}_{\alpha})
\end{align*}
where $T$ denotes tangential derivatives of the given quantities. Thus, one recovers $l^{\alpha \beta}$ and $\partial_n \tilde{A}_{\alpha}$ on $\{x_n=0\}$. By induction, we prove that for $j \geq 1$ 
\begin{multline*}
b_{-j} = -\Big( -\frac{1}{2b_1} \Big)^j ( (\partial_n^{j-1} l^{\alpha \beta}) \omega_{\alpha} \omega_{\beta} + g^{\alpha \beta} (\partial_n^j \tilde{A}_{\alpha}) \omega_{\beta} ) \\
 + T(g^{\alpha \beta}, \partial_n g^{\alpha \beta}, l^{\alpha \beta}, \ldots, \partial_n^{j-2} l^{\alpha \beta}, \tilde{A}_{\alpha}, \ldots, \partial_n^{j-1} \tilde{A}_{\alpha} ).
\end{multline*}
Indeed, this is true for $j = 1$, and assuming this for $j$ one gets 
\begin{align*}
b_{-j-1} &= -\frac{1}{2b_1} \partial_n b_{-j} + T(g^{\alpha \beta}, \partial_n g^{\alpha \beta}, l^{\alpha \beta}, \ldots, \partial_n^{j-1} l^{\alpha \beta}, \tilde{A}_{\alpha}, \ldots, \partial_n^j \tilde{A}_{\alpha} ) \\
 &= -\Big( -\frac{1}{2b_1} \Big)^{j+1} ( (\partial_n^j l^{\alpha \beta}) \omega_{\alpha} \omega_{\beta} + g^{\alpha \beta} (\partial_n^{j+1} \tilde{A}_{\alpha}) \omega_{\beta} ) + T(\,\cdot\,).
\end{align*}
Thus one recovers $\partial_n^j l^{\alpha \beta}$ and $\partial_n^j \tilde{A}_{\alpha}$ on $\{x_n=0\}$ for all $j \geq 0$. The result follows since $\partial^K \tilde{A}_{\alpha} = \partial^K A_{\alpha}$ when $x_n = 0$.
\end{proof}

We may now prove the main result on boundary determination.

\begin{proof}[Proof of Theorem \ref{Bdet:ThmBDet}] Let $\Lambda_{g,A,q} = \Lambda_{\tilde{g},\tilde{A},\tilde{q}}$. Replacing both sets of coefficients by gauge equivalent ones as discussed after Lemma \ref{lemma:boundary_normalization}, we may assume that \eqref{a_normalization} and \eqref{g_normalization} are valid. Then Lemma \ref{lemma:boundary_quantities} implies that the quantities \eqref{boundary_quantities} with and without tildes coincide on $\{x_n = 0\}$. We prove by induction that for $j \geq 0$, one has on $x_n = 0$ 
\begin{equation} \label{boundary_induction_hypothesis}
\partial_n^j q = \partial_n^j \tilde{q}, \quad \partial_n^{j+2} g^{\alpha \beta} = \partial_n^{j+2} \tilde{g}^{\alpha \beta}, \quad g_{\alpha \beta} \partial_n^{j+3} g^{\alpha \beta} = \tilde{g}_{\alpha \beta} \partial_n^{j+3} \tilde{g}^{\alpha \beta}.
\end{equation}

We first note that 
\begin{equation} \label{gdng_two}
g_{\alpha \beta} \partial_n^2 g^{\alpha \beta} = \tilde{g}_{\alpha \beta} \partial_n^2 \tilde{g}^{\alpha \beta} \quad \text{on } x_n = 0.
\end{equation}
This follows from \eqref{g_normalization} for $g$ and $\tilde{g}$, since $\partial_n^j g^{\alpha \beta} = \partial_n^j \tilde{g}^{\alpha \beta}$ on $x_n = 0$ for $j = 0,1$. Note also that \eqref{g_normalization} implies 
\begin{equation*}
\partial_n^j k^{\alpha \beta} = \partial_n^{j+1} g^{\alpha \beta} - (g_{\gamma \delta} \partial_n g^{\gamma \delta}) \partial_n^j g^{\alpha \beta} \quad \text{on } x_n = 0,
\end{equation*}
and therefore 
\begin{equation*}
\partial_n^j l^{\alpha \beta} = \frac{1}{4}(\partial_n^{j+2} g^{\alpha \beta} - (g_{\gamma \delta} \partial_n g^{\gamma \delta}) \partial_n^{j+1} g^{\alpha \beta}) + \partial_n^j(q g^{\alpha \beta}) \quad \text{on } x_n = 0.
\end{equation*}

To prove \eqref{boundary_induction_hypothesis} for $j = 0$, we use that $l^{\alpha \beta} = \tilde{l}^{\alpha \beta}$ on $x_n = 0$. This implies, upon multiplying by $g_{\alpha \beta}$ and summing, that
$q = \tilde{q}$ on $x_n = 0$. Here we used \eqref{gdng_two}. Then $l^{\alpha \beta} = \tilde{l}^{\alpha \beta}$ also implies $\partial_n^2 g^{\alpha \beta} = \partial_n^2 \tilde{g}^{\alpha \beta}$ on $x_n = 0$. The equality $g_{\alpha \beta} \partial_n^3 g^{\alpha \beta} = \tilde{g}_{\alpha \beta} \partial_n^3 \tilde{g}^{\alpha \beta}$ follows by using \eqref{g_normalization}.

Assume now that \eqref{boundary_induction_hypothesis} holds for $j \leq k$. Moving to $k+1$, the equality $\partial_n^{k+1} l^{\alpha \beta} = \partial_n^{k+1} \tilde{l}^{\alpha \beta}$ on $x_n = 0$ implies upon multiplying by $g_{\alpha \beta}$, summing, and using the induction hypothesis, that 
 $$g_{\alpha \beta} \partial_n^{k+1} (q g^{\alpha \beta}) = \tilde{g}_{\alpha \beta} \partial_n^{k+1} (\tilde{q} \tilde{g}^{\alpha \beta}) \ \ \ \text{on } x_n = 0.$$
The induction hypothesis again gives $\partial_n^{k+1} q = \partial_n^{k+1} \tilde{q}$ on $x_n = 0$, and $\partial_n^{k+3} g^{\alpha \beta} = \partial_n^{k+3} \tilde{g}^{\alpha \beta}$ then follows by going back to the equality $\partial_n^{k+1} l^{\alpha \beta} = \partial_n^{k+1} \tilde{l}^{\alpha \beta}$ on $x_n = 0$. The last statement in \eqref{boundary_induction_hypothesis} for $j = k+1$ is a consequence of \eqref{g_normalization}. This ends the induction.

The outcome of the above argument is that $g^{\alpha \beta} \simeq \tilde{g}^{\alpha \beta}$, $A_{\alpha} \simeq \tilde{A}_{\alpha}$, and $q \simeq \tilde{q}$ at $p$. This shows the first statement in Theorem \ref{Bdet:ThmBDet}. If $g \simeq \tilde{g}$ at each $p$, then it is easy to obtain $q \simeq \tilde{q}$ at each $p$ from $l^{\alpha \beta} \simeq \tilde{l}^{\alpha \beta}$. Also, the function $\psi$ constructed locally in Lemma \ref{lemma:boundary_normalization} can be obtained globally on $\partial M$ by a suitable partition of unity. Therefore $A \simeq \tilde{A}$ on all of $\partial M$.
\end{proof}
\end{section}
%
%
\appendix
\begin{section}{Riemannian geometry}
In this appendix we include basic definitions and facts which are used throughout the text. For more details see \cite{Jost}.
We are using Einstein's summation convention: repeated indices in lower and upper position are
summed. In the following $(M,g)$ is a Riemannian manifold. When no confusion may occur, we will use the 
following standard notations for the inner product and the norm:
   $$ \langle X,Y \rangle = g(X,Y), \quad |X|=\sqrt{g(X,X)}. $$ 
\begin{subsection}{Connection and Hessian}
The Riemannian metric $g$ on $M$ induces a natural isomorphism between the tangent 
and cotangent bundles given by
\begin{align*}
   \begin{aligned}
      T(M)  &\to T^*(M) \\
      (x,X) &\mapsto (x,X^{\flat}) 
   \end{aligned}
\end{align*}
where $X^{\flat}(Y)=\langle X,Y\rangle$, and whose inverse is
\begin{align*} 
   \begin{aligned}
      T^*(M)  &\to T(M) \\
      (x,\xi) &\mapsto (x,\xi^{\sharp}) 
   \end{aligned}
\end{align*}
where $\xi^{\sharp}$ is defined by $\xi(X)=\langle \xi^{\sharp},X \rangle$. In local coordinates, if the metric is given by
   $$  g = g_{jk} \,dx^{j} \otimes dx^{k}, $$
this reads
\begin{align*}
   X^{\flat} = g_{j k} X^{j} \,dx^{k}, \quad  \xi^{\sharp} = g^{jk} \xi_{j} \d_{x_k}.
\end{align*}
In particular, the gradient field is defined by $\nabla \phi=d \phi^{\sharp}$. The musical isomorphisms allow to lift
the metric to the cotangent bundle. The cotangent bundle is hence naturally endowed with the Riemannian metric 
$g^{-1}$ given in local coordinates by
\begin{align*}
   g^{-1} = g^{jk} \,d\xi_{j} \otimes d\xi_{k}.
\end{align*} 
It is natural to use $\langle \cdot, \cdot\rangle$ and $|\cdot|$ to denote the inner product and the norm
both on the tangent and cotangent bundles.

We denote by $D$ the Levi-Civita connection on $(M,g)$. A connection is a bilinear map on the vector space
of vector fields which satisfies the following conditions:
\begin{enumerate}
   \item[(\textit{i})] $D_{f X} Y = fD_X Y,$ and $D_X(fY) = (L_Xf)Y+fD_X Y $ if $f$ is a smooth function on $M$, 
   \item[(\textit{ii})] $D_XY-D_YX=[X,Y]$. 
\end{enumerate}  
Here $L_X$ is the Lie derivative. On a Riemannian manifold, there is precisely one connection, called the Levi-Civita connection, which is consistent 
with the metric, i.e.~which satisfies
\begin{enumerate}
   \item[(\textit{iii})] $L_X \langle Y,Z\rangle = \langle D_XY,Z \rangle + \langle Y,D_XZ \rangle$.
\end{enumerate}
This connection is determined in local coordinates by
    $$ D_{\d_{x_j}} \d_{x_k} = \Gamma^l_{jk} \d_{x_l} $$
where the Christoffel symbols $\Gamma^{l}_{jk}$ are given by
   \begin{align*}
      \Gamma^{l}_{jk} = \frac{1}{2} g^{lm} \big(\d_{x_{j}}g_{km}+\d_{x_{k}}g_{jm}-\d_{x_{m}}g_{jk}\big).
   \end{align*}
Note that $\Gamma_{jk}^l=\Gamma_{kj}^l$. If $X$ is a vector field on $M$, then the endomorphism
$D_X$ has a unique extension as an endomorphism on the space of tensors satisfying the following requirements:
\begin{enumerate}
   \item[(\textit{i})] $D_X$ is type preserving,
   \item[(\textit{ii})] $D_X (c(S)) = c(D_X S)$ for any tensor $S$ and any contraction $c$, 
   \item[(\textit{iii})] $D_X(S \otimes T) = D_X T \otimes S + T \otimes D_X S$ for any tensors $S, T$. 
\end{enumerate}  
In particular, if $f$ is a function we have $D_X\phi=d\phi(X)=L_X\phi$, and for 1-forms the connection is given
in local coordinates by
   $$ D_{\d_{x_j}} dx^k = -\Gamma^k_{jl} dx^l. $$
The total derivative $DS$ of $S$ is the tensor $DS(X,\,\cdot\,) = D_X S(\,\cdot\,)$.

The Hessian of a smooth function $\phi$ is the symmetric $(2,0)$-tensor $D^2\phi=Dd\phi$. 
The expression of the Hessian in local coordinates is
   \begin{align*}
      D^2 \phi = \Big(\d^2_{x_j x_k}\phi-\Gamma^{l}_{jk} \d_{x_l}\phi\Big) dx^{j} \otimes dx^{k}.
   \end{align*}
The following identities will be useful:
\begin{gather}
   D^{2}\phi(X,Y)= \frac{1}{2} L_{\nabla \phi}g(X,Y)=\langle D_{X}\nabla \phi,Y\rangle, \label{appendix:HessEq} \\
   D^2 \phi(X,X) = \frac{d^2}{dt^2} \varphi(\gamma(t)) \Big|_{t=0}. \label{appendix:HessGeodesic}
\end{gather}
Here $\gamma$ is the geodesic with $\dot{\gamma}(0) = X$.
\end{subsection}
\begin{subsection}{Parallel and Killing fields}
First we recall the following identities for the Lie derivative: if $f$ is a function
and $X$ a vector field then
\begin{align}
\label{appendix:LieForm}
   L_{f X}g&=fL_Xg+ df \otimes X^{\flat} + X^{\flat} \otimes df,
\end{align}
and if $S$ is a $(2,0)$-tensor then 
\begin{align}
\label{appendix:Lie}
   (L_XS)(Y,Z) = L_X(S(Y,Z))-S([X,Y],Z)-S(Y,[X,Z]). 
\end{align}
\begin{deft}
   A vector field $X$ in $(M,g)$ is a Killing field if
       $$ L_X g = 0. $$ 
\end{deft}
\noindent Note that \eqref{appendix:Lie} implies 
\begin{equation*}
   (L_Xg)(Y,Z) = \langle D_Y X,Z \rangle + \langle Y,D_ZX \rangle. 
\end{equation*}
\begin{deft}
   A vector field $X$ is said to be parallel if its covariant derivative vanishes identically, that is $DX=0$. 
\end{deft}

The following characterization is used in the proof of Theorem \ref{Intro:CharLCW}.

\begin{lem}
\label{appendix:KillLem}
   Let $(M,g)$ be a simply connected Riemannian manifold. A vector field on $(M,g)$ is parallel 
   if and only if it is both a gradient field and a Killing field. 
\end{lem}
\begin{proof}
   In a simply connected manifold, a vector field $X$ is a gradient field if and only if the one form $\omega=X^{\flat}$
   is closed. We have
   \begin{align*}
      d\omega(Y,Z) &= L_Y\big(\omega(Z)\big) -L_Z\big(\omega(Y)\big) -\omega([Y,Z]) \\
      &= L_Y\langle X,Z\rangle -L_Z\langle X,Y \rangle -\langle X,[Y,Z] \rangle \\
      &= \langle D_Y X,Z \rangle+\langle X,D_Y Z \rangle-\langle D_Z X,Y\rangle-\langle X,D_ZY \rangle \\
      &\quad -\langle X,D_Y Z\rangle + \langle X,D_Z Y\rangle. 
   \end{align*}
   Thus, $X$ is a gradient field if and only if for all vector fields $Y,Z$
   \begin{align*}
      \langle D_YX,Z\rangle-\langle D_ZX,Y \rangle=0.
   \end{align*}
   On the other hand, $X$ is a Killing field if and only if
   \begin{align*}
      L_Xg (Y,Z) = \langle D_YX,Z \rangle+ \langle Y,D_ZX \rangle=0
   \end{align*}
   for all vector fields $Y,Z$. The result ensues.
\end{proof}

The next result states that the existence of a unit parallel vector field implies a local product structure on the manifold.

\begin{lem} \label{appendix:localLem}
Let $X$ be a unit parallel vector field in a manifold $(M,g)$. Near any point of $M$, there are local coordinates $x$ such that $X = \partial/\partial x_1$ and the metric has the form 
\begin{equation*}
g(x_1,x') = \left( \begin{array}{cc} 1 & 0 \\ 0 & g_0(x') \end{array} \right).
\end{equation*}
Conversely, if such coordinates exist then $X = \partial/\partial x_1$ is unit parallel.
\end{lem}
\begin{proof}
Let $X$ be unit parallel and let $\Gamma$ be the distribution orthogonal to $X$. If $V, W$ are vector fields orthogonal to $X$ then 
\begin{equation*}
\langle [V,W], X \rangle = \langle D_V W - D_W V, X \rangle = V \langle W,X \rangle - W \langle V,X \rangle = 0.
\end{equation*}
Thus $\Gamma$ is involutive, and by the Frobenius theorem there is a hypersurface $S$ (through any point of $M$) which is normal to $X$. Let $x' \mapsto q(x')$ be local coordinates on $S$, and let $(x_1,x') \mapsto \exp_{q(x')}(x_1 X(q(x'))$ be corresponding semigeodesic coordinates. In fact integral curves of $X$ are geodesics (if $\dot{\gamma}(t) = X(\gamma(t))$ then $D_{\dot{\gamma}} \dot{\gamma} = 0$), so $X = \partial_1$. Then 
\begin{equation*}
g(x_1,x') = \left( \begin{array}{cc} 1 & 0 \\ 0 & g_0(x_1,x') \end{array} \right).
\end{equation*}
If $j,k \geq 2$ then $\partial_1 g_{jk} = \langle D_{\partial_1} \partial_j, \partial_k \rangle + \langle \partial_j, D_{\partial_1} \partial_k \rangle = 0$ since $\partial_1$ is parallel. Therefore $g_0 = g_0(x')$.

The converse follows since $D_{\partial_j} \partial_1 = 0$ by a direct computation.
\end{proof}

Finally, we need the definition of conformal Killing fields.
\begin{deft}
   A vector field $X$ in $(M,g)$ is called a conformal Killing field if
       $$ L_X g = \lambda g. $$
\end{deft}
Note that if $L_X g = \lambda g$, by taking traces one has $\lambda=\frac{2}{n} \mathop{\rm div} X$. The notion of conformal Killing field is of course invariant under conformal change of the metric. 

\end{subsection}

\begin{subsection}{Submanifolds}

Let $(M,g)$ be a Riemannian manifold and let $S$ be an embedded hypersurface in $M$. Fix a unit vector field $N$ which is normal to $S$. The second fundamental form of $S$ is defined by 
\begin{equation*}
\ell(X,Y) = \langle D_X N, Y \rangle,
\end{equation*}
where $X$ and $Y$ are vector fields tangent to $S$, and $D$ is the Levi-Civita connection in $(M,g)$. The eigenvalues of the symmetric bilinear form $\ell$ are the principal curvatures of $S$; their sign depends on the choice of normal.

\begin{deft}
\label{appendix:umbilical}
A point of a hypersurface is called umbilical if all the principal curvatures are equal at that point. A hypersurface is called totally umbilical if every point is umbilical.
\end{deft}

\begin{deft}
\label{appendix:strictlyconvex}
A hypersurface is called strictly convex if the second fundamental form is positive definite.
\end{deft}

\end{subsection}

\begin{subsection}{Curvature tensors}

Next we consider curvature tensors on a Riemannian manifold $(M,g)$. The Riemann curvature endomorphism is a $(3,1)$-tensor on $M$, defined by 
\begin{equation*}
R(X,Y)Z = D_X D_Y Z - D_Y D_X Z - D_{[X,Y]} Z,
\end{equation*}
whenever $X,Y,Z$ are vector fields on $M$. By lowering indices, one obtains the Riemann curvature tensor which is the $(4,0)$-tensor 
\begin{equation*}
R(X,Y,Z,W) = \langle R(X,Y)Z, W \rangle.
\end{equation*}
In local coordinates, with coordinate vector fields $\partial_a = \partial/\partial x_a$ and with $D_a = D_{\partial_a}$, the components of the curvature tensor are given by 
\begin{equation*}
R_{abcd} = \langle (D_a D_b - D_b D_a) \partial_c, \partial_d \rangle.
\end{equation*}
By taking traces of the Riemann curvature tensor, we obtain the Ricci tensor which is a symmetric $(2,0)$-tensor whose components are 
\begin{equation*}
R_{bc} = g^{ad} R_{abcd}.
\end{equation*}
The scalar curvature is the function 
\begin{equation*}
\mathrm{Scal} = g^{bc} R_{bc}.
\end{equation*}
A manifold $(M,g)$ is said to be flat if the Riemann curvature tensor vanishes identically. Euclidean space is flat, and any flat manifold is locally isometric to a subset of Euclidean space.

Finally, we wish to define the conformal curvature tensors. First consider the rho-tensor, which is a symmetric $2$-tensor given in components by 
\begin{equation*}
P_{ab} = \frac{1}{n-2} \Big( R_{ab} - \frac{\mathrm{Scal}}{2(n-1)} g_{ab} \Big).
\end{equation*}
The Weyl tensor of $(M,g)$ is the $4$-tensor 
\begin{equation*}
W_{abcd} = R_{abcd} + P_{ac} g_{bd} + P_{bd} g_{ac} - P_{bc} g_{ad} - P_{ad} g_{bc},
\end{equation*}
and the Cotton tensor is the $3$-tensor 
\begin{equation*}
C_{abc} = D_a P_{bc} - D_b P_{ac}.
\end{equation*}
If the metric $g$ is replaced by a conformal metric $cg$, then the Weyl tensor transforms as $W_{cg} = c W_g$. If $n = 3$ then $W \equiv 0$, but one has $C_{cg} = C_g$.

A manifold $(M,g)$ is called conformally flat if some conformal manifold $(M,cg)$ is flat. Any $2$-dimensional manifold is locally conformally flat. A $3$-dimensional manifold is locally conformally flat if and only if its Cotton tensor vanishes identically, and a manifold of dimension $n \geq 4$ is locally conformally flat if and only if the Weyl tensor vanishes (see \cite{eisenhart}).

\end{subsection}
\end{section}
%
%

%

\begin{thebibliography}{99} %
%
\bibitem{An} Yu.~E.~Anikonov, \textit{Some methods for the study of multidimensional inverse problems for differential equations}, Nauka Sibirsk. Otdel, Novosibirsk (1978).
\bibitem{ALP} K.~Astala, M.~Lassas, L.~P\"aiv\"arinta, \textit{Calder{\'o}n's inverse problem for anisotropic conductivity in the plane}, Comm. Partial Differential Equations, \textbf{30} (2005), 207--224.
\bibitem{AP} K.~Astala, L.~P\"aiv\"arinta, \textit{Calder{\'o}n's inverse conductivity problem in the plane}, Ann. of Math., \textbf{163} (2006), 265--299.
\bibitem{BB} D.~C.~Barber, B.~H.~Brown, \textit{Progress in electrical impedance tomography}, in Inverse problems in partial differential equations, edited by D.~Colton, R.~Ewing, and W.~Rundell, SIAM, Philadelphia (1990), 151--164.
\bibitem{BrU} R.~M.~Brown, G.~Uhlmann, \textit{Uniqueness in the inverse conductivity problem for nonsmooth conductivities in two dimensions}, Comm. Partial Differential Equations, \textbf{22} (1997), 1009--1027.
\bibitem{DS} M.~Dimassi, J.~Sj\"ostrand, \textit{Spectral asymptotics in the semi-classical limit}, Cambridge University Press, 1999.
\bibitem{DSFKSU} D.~Dos Santos Ferreira, C.~E.~Kenig, J.~Sj\"ostrand, G.~Uhlmann,  \textit{Determining a magnetic Schr\"odinger operator
   from partial Cauchy data}, Comm. Math. Phys., \textbf{271} (2007), 467--488.
\bibitem{eisenhart} L.~Eisenhart, \emph{{Riemannian geometry}}, 2nd printing, Princeton University Press, 1949.
\bibitem{GS} C.~Guillarmou, A.~Sa Barreto, \textit{Inverse problems for Einstein manifolds}, preprint (2007), arXiv:0710.1136.
\bibitem{H3} L.~H\"ormander, \textit{The Analysis of Linear Partial Differential Operators III}, Springer-Verlag, 1985.
\bibitem{I} H.~Isozaki, \textit{Inverse spectral problems on hyperbolic manifolds and their applications to inverse boundary value problems in {E}uclidean space}, Amer. J. Math., \textbf{126} (2004), 1261--1313.
\bibitem{Jost} J.~Jost, \textit{Riemannian geometry and geometric analysis}, Springer-Verlag, 2002.
\bibitem{KSU} C.~E.~Kenig, J.~Sj\"ostrand, G.~Uhlmann,  \textit{The Calder\'on problem with partial data}, Ann. of Math., \textbf{165} (2007), 567--591.
\bibitem{KnSa} K.~Knudsen, M.~Salo, \textit{Determining non-smooth first order terms from partial boundary measurements}, Inverse Problems and Imaging,
   \textbf{1} (2007), 349--369.
\bibitem{KV} R.~Kohn, M.~Vogelius, \textit{Identification of an unknown conductivity by means of measurements at the boundary}, in Inverse Problems, edited by D.~McLaughlin, SIAM-AMS Proc. No. 14, Amer. Math. Soc., Providence (1984), 113--123.
\bibitem{LTU} M.~Lassas, M.~Taylor, G.~Uhlmann, \textit{The Dirichlet-to-Neumann map for complete Riemannian manifolds with boundary}, Comm. Anal. Geom., \textbf{11} (2003), 207--221.
\bibitem{LaU} M.~Lassas, G.~Uhlmann, \textit{On determining a Riemannian manifold from the Dirichlet-to-Neumann map}, Ann. Sci. {\'E}cole Norm. Sup., \textbf{34} 
   (2001), 771--787.
\bibitem{LeU} J.~Lee, G.~Uhlmann, \textit{Determining anisotropic real-analytic conductivities by boundary measurements}, Comm. Pure Appl. Math., 
   \textbf{42} (1989), 1097--1112.
\bibitem{Li} W.~Lionheart, \textit{Conformal uniqueness results in anisotropic electrical impedance imaging}, Inverse Problems, \textbf{13} (1997), 125–-134.
\bibitem{M} R.~G.~Mukhometov, \textit{The reconstruction problem of a two-dimensional Riemannian metric, and integral geometry} (Russian), Dokl. Akad. Nauk SSSR, \textbf{232} (1977), 32–-35.
\bibitem{N} A.~Nachman, \textit{Global uniqueness for a two-dimensional inverse boundary value problem}, Ann. of Math., \textbf{143} (1996), 71--96.
\bibitem{NSU} G.~Nakamura, Z.~Sun, G.~Uhlmann, \textit{Global identifiability for an inverse problem for the Schr\"odinger equation 
   in a magnetic field}, Math. Ann., \textbf{303} (1995), 377--388.
\bibitem{PW} L.~E.~Payne, H.~F.~Weinberger, \textit{An optimal Poincar{\'e} inequality for convex domains}, Arch. Rat. Mech. Anal., \textbf{5} (1960), 286--292.
\bibitem{Pe} P.~Petersen, \textit{Riemannian geometry}, Springer-Verlag, 1998.
\bibitem{Sa} M.~Salo, \textit{Inverse boundary value problems for the magnetic Schr{\"o}dinger equation}, J. Phys. Conf. Series, \textbf{73} (2007), 012020.
\bibitem{SaTz} M.~Salo, L.~Tzou, \textit{Carleman estimates and inverse problems for Dirac operators}, preprint (2007), arXiv:0709.2282.
\bibitem{Sh} V.~Sharafutdinov, \textit{Integral geometry of tensor fields}, in \textit{Inverse and Ill-Posed Problems Series}, VSP, Utrecht, 1994.
\bibitem{Sh2} V.~Sharafutdinov, \textit{On emission tomography of inhomogeneous media}, SIAM J. Appl. Math., \textbf{55} (1995), 707--718.
\bibitem{SuU} Z.~Sun, G.~Uhlmann, \textit{Generic uniqueness for an inverse boundary value problem}, Duke Math. J., \textbf{62} (1991), 131--155.
\bibitem{SuU2} Z.~Sun, G.~Uhlmann, \textit{Anisotropic inverse problems in two dimensions}, Inverse Problems, \textbf{19} (2003), 1001--1010.
\bibitem{Sy} J.~Sylvester, \textit{An anisotropic inverse boundary value problem}, Comm. Pure Appl. Math., \textbf{43} (1990), 201--232.
\bibitem{SU} J.~Sylvester, G.~Uhlmann, \textit{A global uniqueness theorem for an inverse boundary value problem}, Ann. of Math., \textbf{125} 
   (1987), 153--169.
\bibitem{SU2} J.~Sylvester, G.~Uhlmann, \textit{Inverse boundary value problems at the boundary -- continuous dependence}, Comm. Pure Appl. Math., 
   \textbf{41} (1988), 197--219.
\bibitem{T} M.~Taylor, \textit{Pseudodifferential operators}, Princeton University Press, 1981.
\end{thebibliography}
\end{document}